\RequirePackage{ifpdf}
\ifpdf 
\documentclass[pdftex]{sigma}
\else
\documentclass{sigma}
\fi

\numberwithin{equation}{section}

\begin{document}
\allowdisplaybreaks

\renewcommand{\PaperNumber}{052}

\FirstPageHeading

\ShortArticleName{Polynomials Associated with Dihedral Groups}

\ArticleName{Polynomials Associated with Dihedral Groups}

\Author{Charles F. DUNKL}

\AuthorNameForHeading{C.F. Dunkl}

\Address{Department of Mathematics, University of Virginia,
Charlottesville, VA 22904-4137, USA}
\Email{\href{mailto:cfd5z@virginia.edu}{cfd5z@virginia.edu}}
\URLaddress{\url{http://www.people.virginia.edu/~cfd5z}}

\ArticleDates{Received February 06, 2007; Published online March 22, 2007}

\Abstract{There is a commutative algebra of dif\/ferential-dif\/ference operators, with two
parameters, associated to any dihedral group with an even number of
ref\/lections. The intertwining operator relates this algebra to the algebra of
partial derivatives. This paper presents an explicit form of the action of the
intertwining operator on polynomials by use of harmonic and Jacobi
polynomials. The last section of the paper deals with parameter values for
which the formulae have singularities.}


\Classification{33C45; 33C80; 20F55}

\section{Introduction}

The dihedral group of type $I_{2}\left(  2s\right)  $ acts on $\mathbb{R}^{2}$, contains $2s$ ref\/lections and the rotations through angles of
$\frac{m\pi}{s}$ for $1\leq m\leq2s-1$, and is of order $4s$, where $s$ is a
positive integer. It is the symmetry group of the regular $2s$-gon and has two
conjugacy classes of ref\/lections (the mirrors passing through midpoints of
pairs of opposite edges and those joining opposite vertices). There is an
associated commutative algebra of dif\/ferential-dif\/ference (``\textit{Dunkl}'') operators with two parameters, denoted by
$\kappa_{0}$, $\kappa_{1}$. It is convenient to use complex coordinates for
$\mathbb{R}^{2}$, that is, $z=x_{1}+\mathrm{i}x_{2}$, $\overline{z}=x_{1}-\mathrm{i}x_{2}$.
Notations like $f\left(  z\right)  $ will be understood as functions of
$z$, $\overline{z}$; except that $f\left(  \overline{z},z\right)  $ will be used
to indicate the result of interchanging $z$ and $\overline{z}$. Let $\mathbb{N}$,
$\mathbb{N}_{0}$, $\mathbb{Q}$ denote the sets of positive integers, nonnegative integers and rational
numbers, respectively. Let $\omega=e^{\mathrm{i}\pi/s}$, then the ref\/lections
in the group are $\left(  z,\overline{z}\right)  \mapsto\left(  \overline
{z}\omega^{m},z\omega^{-m}\right)$, $0\leq m<2s$ and the rotations are $\left(
z,\overline{z}\right)  \mapsto\left(  z\omega^{m},\overline{z}\omega
^{-m},\right)$, $1\leq m<2s$. Note that $f\left(  \overline{z}\omega
^{m}\right)  $ is the abbreviated form of $f\left(  \overline{z}\omega
^{m},z\omega^{-m}\right)  $. The dif\/ferential-dif\/ference operators are def\/ined
by
\begin{gather*}
Tf\left(  z\right)     :=\frac{\partial}{\partial z}f\left(  z\right)
+\kappa_{0}\sum_{j=0}^{s-1}\frac{f\left(  z\right)  -f\left(  \overline
{z}\omega^{2j}\right)  }{z-\overline{z}\omega^{2j}}+\kappa_{1}\sum_{j=0}%
^{s-1}\frac{f\left(  z\right)  -f\left(  \overline{z}\omega^{2j+1}\right)
}{z-\overline{z}\omega^{2j+1}},\\
\overline{T}f\left(  z\right)     :=\frac{\partial}{\partial\overline{z}%
}f\left(  z\right)  -\kappa_{0}\sum_{j=0}^{s-1}\frac{f\left(  z\right)
-f\left(  \overline{z}\omega^{2j}\right)  }{z-\overline{z}\omega^{2j}}%
\omega^{2j}-\kappa_{1}\sum_{j=0}^{s-1}\frac{f\left(  z\right)  -f\left(
\overline{z}\omega^{2j+1}\right)  }{z-\overline{z}\omega^{2j+1}}\,\omega^{2j+1},
\end{gather*}
for polynomials $f\left(  z\right)  $. (The second formula implicitly uses the
relation $-\frac{\omega^{m}}{z-\overline{z}\omega^{m}}=\frac{1}{\overline
{z}-z\omega^{-m}}$.) The key fact is that $T$ and $\overline{T}$ commute. The
explicit action of $T$ and $\overline{T}$ on monomials is given by
\begin{gather}
Tz^{a}\overline{z}^{b}    =az^{a-1}\overline{z}^{b}+s\sum_{j=0}^{\left\lfloor
\left(  a-b-1\right)  /s\right\rfloor }\big(  \kappa_{0}+\left(  -1\right)
^{j}\kappa_{1}\big)  z^{a-1-js}\overline{z}^{b+js},\label{TZ}\\
\overline{T}z^{a}\overline{z}^{b}    =bz^{a}\overline{z}^{b-1}-s\sum
_{j=1}^{\left\lfloor \left(  a-b\right)  /s\right\rfloor }\big(  \kappa
_{0}+\left(  -1\right)  ^{j}\kappa_{1}\big)  z^{a-js}\overline{z}^{b-1+js},
\label{TZB}%
\end{gather}
for $a\geq b$; the relations remain valid when both $\left(  z,\overline
{z}\right)  $ and $\left(  T,\overline{T}\right)  $ are interchanged. The
Laplacian is $4T\overline{T}$. These results are from \cite[Section~3]{D1}. The
harmonic polynomials and formulae (\ref{TZ}) and (\ref{TZB}) also appear in
Berenstein and Burman \cite[Section~2]{BB}. The aim of this paper is to f\/ind an
explicit form of the intertwining operator $V$. This is the unique linear
transformation that maps homogeneous polynomials to homogeneous polynomials of
the same degree and satisf\/ies
\[
TVf\left(  z\right)  =V\frac{\partial}{\partial z}f\left(  z\right)
,\qquad \overline{T}Vf\left(  z\right)  =V\frac{\partial}{\partial\overline{z}%
}f\left(  z\right)  ,\qquad V1=1.
\]
The operator was def\/ined for general f\/inite ref\/lection groups in \cite{D3}.
R\"{o}sler \cite{R} proved that $V$ is a~positive operator when $\kappa
_{0},\kappa_{1}>0$; this roughly means that if a polynomial $f$ satisf\/ies
$f\left(  y\right)  \geq0$ for all $y$ with $\left\Vert y\right\Vert <R$ (for
some $R$) then $Vf\left(  y\right)  \geq0$ on the same set. The present paper
does not shed light on the positivity question since the formulae are purely
algebraic. In Section 5 the special case $-\left(  \kappa_{0}+\kappa
_{1}\right)  \in%
\mathbb{N}
$ is considered in more detail. These values of $\left(  \kappa_{0},\kappa
_{1}\right)  $ are apparent singularities in the expressions for
$Vz^{a}\overline{z}^{b}$ which are found in Section 4. The book by Y. Xu and
the author \cite{DX} is a convenient reference for the background of this paper.

In a way, to f\/ind $Vz^{a}\overline{z}^{b}$ only requires to solve a set of
equations involving $Vz^{j}\overline{z}^{k}$ for $0\leq j\leq a$, $0\leq k\leq
b$. This can be implemented in computer algebra for small $a$, $b$ but it is not
really an explicit description. For example, by direct computation we f\/ind
that%
\begin{gather*}
Vz^{2}    =\frac{\left(  \kappa_{0}+\kappa_{1}+1\right)  z^{2}+\left(
\kappa_{0}-\kappa_{1}\right)  \overline{z}^{2}}{\left(  2\kappa_{0}+1\right)
\left(  2\kappa_{1}+1\right)  \left(  2\kappa_{0}+2\kappa_{1}+1\right)
}, \qquad s=2,\\
Vz^{2}    =\dfrac{2z^{2}}{\left(  s\kappa_{0}+s\kappa_{1}+1\right)  \left(
s\kappa_{0}+s\kappa_{1}+2\right)  }, \qquad s>2.
\end{gather*}
The idea is to f\/ind the harmonic expansion of $Vz^{a}\overline{z}^{b}$;
suppose $f\left(  z\right)  $ is (real-) homogeneous of degree $n$ then there
is a unique expansion $f\left(  z\right)  =\sum\limits_{j=0}^{\left\lfloor
n/2\right\rfloor }\left(  z\overline{z}\right)  ^{j}f_{n-2j}\left(  z\right)
$ where $f_{n-2j}$ is homogeneous of degree $n-2j$ and is harmonic, that is,
$T\overline{T}f_{n-2j}=0$, for $0\leq j\leq n/2$. There is some more
information easily available for the expansion of $Vz^{a}\overline{z}^{b}$.
Let $n=a+b$ and suppose $Vz^{a}\overline{z}^{b}=\sum\limits_{j=0}^{n}c_{j}%
z^{n-j}\overline{z}^{j}$ for certain coef\/f\/icients $c_{j}$. Because $V$
commutes with the action of the group we deduce that
\[
V\left(  \left(  \omega z\right)  ^{a}\left(  \overline{\omega z}\right)
^{b}\right)  =\omega^{a-b}\sum_{j=0}^{n}c_{j}z^{n-j}\overline{z}^{j}%
=\sum_{j=0}^{n}\omega^{n-2j}c_{j}z^{n-j}\overline{z}^{j};
\]
thus $c_{j}\neq0$ implies $n-2j\equiv a-b\operatorname{mod}\left(  2s\right)
$ or $j\equiv b\operatorname{mod}s$. Further%
\[
V\left(  \overline{z}^{a}z^{b}\right)  =\sum_{j=0}^{n}c_{j}\overline{z}%
^{n-j}z^{j},
\]
so it will suf\/f\/ice to determine $Vz^{a}\overline{z}^{b}$ for $a\geq b$. We
will use the Poisson kernel to calculate the polynomials denoted $K_{n}\left(
x,y\right)  :=V^{x}\left(  \frac{1}{n!}\left(  x_{1}y_{1}+x_{2}y_{2}\right)
^{n}\right)  $ (see \cite[p.~1219]{D4}), where $y\in\mathbb{R}^{2}$ and $V^{x}$
acts on the variable $x$. Thus $Vx_{1}^{n-j}x_{2}^{j}$ is $j!\left(
n-j\right)  !$ times the coef\/f\/icient of $y_{1}^{n-j}y_{2}^{j}$ in
$K_{n}\left(  x,y\right)  $. This is adapted to complex coordinates by setting
$w=y_{1}+\mathrm{i}y_{2}$, in which case $x_{1}y_{1}+x_{2}y_{2}=\frac{1}%
{2}\left(  z\overline{w}+\overline{z}w\right)  $.

\section{The Poisson kernel}

Actually it is only the series expansion of this kernel that is used. For now
we assume $\kappa_{0},\kappa_{1}\geq0$. The measure on the circle
$\mathbb{T}:=\left\{  e^{\mathrm{i}\theta}:-\pi<\theta\leq\pi\right\}  $
associated to the group $I_{2}\left(  2s\right)  $ and the operators
$T$, $\overline{T}$ is
\[
d\mu\big(  e^{\mathrm{i}\theta}\big)  :=\frac{1}{2B\left(  \kappa_{0}%
+\frac{1}{2},\kappa_{1}+\frac{1}{2}\right)  }\left(  \sin^{2}s\theta\right)
^{\kappa_{0}}\left(  \cos^{2}s\theta\right)  ^{\kappa_{1}} d\theta.
\]
Suppose $g$ is a function of $t=\cos2s\theta$ then
\[
\int_{\mathbb{T}}g\left(  t\left(  \theta\right)  \right)  d\mu\big(
e^{\mathrm{i}\theta}\big)  =\frac{2^{-\kappa_{0}-\kappa_{1}}}{B\left(
\kappa_{0}+\frac{1}{2},\kappa_{1}+\frac{1}{2}\right)  }\int_{-1}^{1}g\left(
t\right)  \left(  1-t\right)  ^{\kappa_{0}-1/2}\left(  1+t\right)
^{\kappa_{1}-1/2}dt.
\]
The inner product in $L^{2}\left(  \mathbb{T},\mu\right)  $ is
\[
\left\langle f,g\right\rangle :=\int_{\mathbb{T}}f\left(  z\right)
\overline{g\left(  z\right)  }d\mu\left(  z\right)
\]
and $\left\Vert f\right\Vert :=\left\langle f,f\right\rangle ^{1/2}$.
Throughout the polynomials under consideration have real coef\/f\/icients so that
$\overline{g\left(  z,\overline{z}\right)  }=g\left(  \overline{z},z\right)
$. By the group invariance of $\mu$ the integral $\int_{\mathbb{T}}%
z^{a}\overline{z}^{b}d\mu\left(  z\right)  $is real-valued when $a\equiv
b\operatorname{mod}\left(  2s\right)  $ and vanishes otherwise. There is an
orthogonal decomposition $L^{2}\left(  \mathbb{T},\mu\right)  =\sum\limits
_{n=0}^{\infty}\oplus\mathcal{H}_{n}$; for $n>0$ each $\mathcal{H}_{n}$ is of
dimension two and consists of the polynomials in $z,\overline{z}$ (real-)
homogeneous of degree $n$ and annihilated by $T\overline{T}$ (the harmonic
property), while $\mathcal{H}_{0}$ consists of the constant functions. The
Poisson kernel is the reproducing kernel for harmonic polynomials (for more
details see \cite{D2,D4}). Xu \cite{X} investigated relationships
between harmonic polynomials, the intertwining operator and the Poisson kernel
for the general ref\/lection group. The paper of Scalas~\cite{S} concerns
boundary value problems for the dihedral groups. The projection of the kernel
onto $\mathcal{H}_{n}$ is denoted by $P_{n}\left(  z,w\right)  $ and
satisf\/ies
\[
\int_{\mathbb{T}}P_{n}\left(  z,w\right)  g\left(  w\right)  d\mu\left(
w\right)  =g\left(  z\right)
\]
for each polynomial $g\in\mathcal{H}_{n}$. There is a formula for $P_{n}$ in
terms of $\left\{  K_{n-2j}:0\leq j\leq\frac{n}{2}\right\}  $ (see
\cite[p.~1224]{D4}) which can be inverted. In the present case
\[
P_{n}\left(  z,w\right)  =\sum_{j=0}^{\left\lfloor n/2\right\rfloor }%
\frac{\left(  \gamma_{0}\right)  _{n}}{j!\left(  2-n-\gamma_{0}\right)  _{j}%
}2^{n-2j}\left(  z\overline{z}w\overline{w}\right)  ^{j}K_{n-2j}\left(
z,w\right)  ,
\]
where $\gamma_{0}=s\kappa_{0}+s\kappa_{1}+1$. The inverse relation is%
\begin{equation}
K_{n}\left(  z,w\right)  =2^{-n}\sum_{j=0}^{\left\lfloor n/2\right\rfloor
}\frac{1}{j!\left(  \gamma_{0}\right)  _{n-j}}\left(  z\overline{z}%
w\overline{w}\right)  ^{j}P_{n-2j}\left(  z,w\right)  . \label{P2K}%
\end{equation}
This is a consequence of the following:

\begin{proposition}
Suppose there are two sequences $\left\{  \xi_{n}:n\in%
\mathbb{N}
_{0}\right\}  $ and $\left\{  \eta_{n}:n\in%
\mathbb{N}
_{0}\right\}  $ in a vector space over $\mathbb{%
\mathbb{Q}
}\left(  \gamma_{0}\right)  $ where $\gamma_{0}$ is transcendental, then
\[
\xi_{n}=\sum_{j=0}^{\left\lfloor n/2\right\rfloor }\frac{\left(  \gamma
_{0}\right)  _{n}}{j!\left(  2-n-\gamma_{0}\right)  _{j}}\eta_{n-2j},\qquad n\in
\mathbb{N}
_{0},
\]
if and only if
\[
\eta_{n}=\sum_{j=0}^{\left\lfloor n/2\right\rfloor }\frac{1}{j!\left(
\gamma_{0}\right)  _{n-j}}\xi_{n-2j},\qquad n\in
\mathbb{N}
_{0}.
\]
\end{proposition}

\begin{proof}
Consider the matrices $A$ and $B$ def\/ined by $\xi_{n}=\sum_{j}A_{jn}\eta
_{j}$, $\eta_{n}=\sum_{j}B_{jn}\xi_{j}$; these matrices are triangular and the
diagonal entries are nonzero, hence they are nonsingular. It suf\/f\/ices to show
$B$ is a one-sided inverse of $A$; this is actually f\/inite-dimensional linear
algebra, since one can truncate to the range $0\leq n,j\leq M$ for any $M\in
\mathbb{N}
$. Indeed
\begin{gather*}
  \sum_{j=0}^{\left\lfloor n/2\right\rfloor }\frac{\left(  \gamma_{0}\right)
_{n}}{j!\left(  2-n-\gamma_{0}\right)  _{j}}\sum_{i=0}^{\left\lfloor
n/2-j\right\rfloor }\frac{1}{i!\left(  \gamma_{0}\right)  _{n-2j-i}}%
\xi_{n-2j-2i}\\
\qquad{}  =\sum_{k=0}^{\left\lfloor n/2\right\rfloor }\xi_{n-2k}\frac{\left(
\gamma_{0}\right)  _{n}}{k!\left(  \gamma_{0}\right)  _{n-k}}\sum_{j=0}%
^{k}\frac{\left(  -k\right)  _{j}\left(  1-n-\gamma_{0}+k\right)  _{j}%
}{j!\left(  2-n-\gamma_{0}\right)  _{j}}\\
\qquad{}  =\sum_{k=0}^{\left\lfloor n/2\right\rfloor }\xi_{n-2k}\frac{\left(
\gamma_{0}\right)  _{n}\left(  1-k\right)  _{k}}{k!\left(  \gamma_{0}\right)
_{n-k}\left(  2-n-\gamma_{0}\right)  _{k}}=\xi_{n};
\end{gather*}
using the substitution $i=k-j$ we obtain $\left(  \gamma_{0}\right)
_{n-2j-i}=\left(  \gamma_{0}\right)  _{n-k-j}=\frac{\left(  -1\right)
^{j}\left(  \gamma_{0}\right)  _{n-k}}{\left(  1-n-\gamma_{0}+k\right)  _{j}}$
and $\frac{1}{i!}=\left(  -1\right)  ^{j}\frac{\left(  -k\right)  _{j}}{k!}$;
the sum over $j$ is found by the Chu--Vandermonde formula.
\end{proof}

Set $\xi_{n}=\frac{P_{n}\left(  z,w\right)  }{\left(  z\overline{z}%
w\overline{w}\right)  ^{n/2}}$ and $\eta_{n}=\frac{2^{n}K_{n}\left(
z,w\right)  }{\left(  z\overline{z}w\overline{w}\right)  ^{n/2}}$ for $n\in%
\mathbb{N}
_{0}$ to prove equation (\ref{P2K}).

Suppose for each $n\in%
\mathbb{N}
$ there exist a basis $\left\{  h_{n1},h_{n2}\right\}  $ and a biorthogonal
basis $\left\{  g_{n1},g_{n2}\right\}  $ for $\mathcal{H}_{n}$ with real
coef\/f\/icients in $z,\overline{z}$ (so $\overline{h_{n1}\left(  z,\overline
{z}\right)  }=h_{n1}\left(  \overline{z},z\right)  $, for example). Thus
$\left\langle h_{ni},g_{nj}\right\rangle =\delta_{ij}/\lambda_{ni},$with
structural constants $\lambda_{ni}$. Then
\begin{equation}
P_{n}\left(  z,w\right)  =\sum_{i=1}^{2}\lambda_{ni}h_{ni}\left(
z,\overline{z}\right)  g_{ni}\left(  \overline{w},w\right)  . \label{pkn}
\end{equation}
Once this is made suf\/f\/iciently explicit we can compute $K_{n}\left(
z,w\right)  $ and $Vz^{n-j}\overline{z}^{j}$. The description of harmonic
polynomials is in terms of the case $s=1$ (corresponding to the group
$I_{2}\left(  2\right)  =\mathbb{Z}_{2}\times\mathbb{Z}_{2}$). In terms of
Jacobi polynomials the polynomials annihilated by $\overline{T}$ are:
\begin{gather}
f_{2n}\left(  re^{\mathrm{i}\theta}\right)     :=r^{2n}P_{n}^{\left(
\kappa_{0}-\frac{1}{2},\kappa_{1}-\frac{1}{2}\right)  }\left(  \cos
2\theta\right)  +\frac{\mathrm{i}}{2}\left(  r^{2}\sin2\theta\right)
r^{2n-2}P_{n-1}^{\left(  \kappa_{0}+\frac{1}{2},\kappa_{1}+\frac{1}{2}\right)
}\left(  \cos2\theta\right)  ,\label{f2nj}\\
f_{2n+1}\left(  re^{\mathrm{i}\theta}\right)    :=\left(  n+\kappa_{0}%
+\frac{1}{2}\right)  r\cos\theta~r^{2n}P_{n}^{\left(  \kappa_{0}-\frac{1}%
{2},\kappa_{1}+\frac{1}{2}\right)  }\left(  \cos2\theta\right)  \label{f2n1j}%
\\
\phantom{f_{2n+1}\left(  re^{\mathrm{i}\theta}\right)    :=}{}
+\mathrm{i}\left(  n+\kappa_{1}+\frac{1}{2}\right)  r\sin\theta~r^{2n}%
P_{n}^{\left(  \kappa_{0}+\frac{1}{2},\kappa_{1}-\frac{1}{2}\right)  }\left(
\cos2\theta\right)  ;\nonumber
\end{gather}
where the subscript indicates the degree of homogeneity, (clearly $f_{n}$ is a
polynomial with real coef\/f\/icients in $z$, $\overline{z}$; $\cos2\theta=\left(
z^{2}+\overline{z}^{2}\right)  /\left(  2z\overline{z}\right)  $ and
$\frac{\mathrm{i}}{2}\left(  r^{2}\sin2\theta\right)  =\frac{1}{4}\left(
z^{2}-\overline{z}^{2}\right)  $). The real and imaginary parts form a basis
for the harmonic polynomials. Specif\/ically let%
\begin{gather*}
f_{n}^{0}\left(  z\right)     :=\operatorname{Re}f_{n}\left(  z\right)  ,\qquad
f_{n}^{1}\left(  z\right)     :=\mathrm{i}\operatorname{Im}f_{n}\left(
z\right)  .
\end{gather*}
This implies that both $f_{n}^{0}$ and $f_{n}^{1}$ have real coef\/f\/icients in
$z,\overline{z}$ and $f_{n}^{0}\left(  \overline{z},z\right)  =f_{n}%
^{0}\left(  z,\overline{z}\right)  ,~f_{n}^{1}\left(  \overline{z},z\right)
=-f_{n}^{1}\left(  z,\overline{z}\right)  $. When $s>1$ and $1\leq t<s$ it is
known \cite[p.~182]{D1} that $\left\{  z^{t}f_{n}\left(  z^{s}\right)
,\overline{z}^{t}f_{n}\left(  \overline{z}^{s}\right)  \right\}  $ is an
orthogonal basis for $\mathcal{H}_{ns+t}$ for $n\geq0$. Henceforth we denote
$h_{ns+t,1}\left(  z\right)  =g_{ns+t,1}\left(  z\right)  =z^{t}f_{n}\left(
z^{s}\right)  =\overline{h_{ns+t,2}}=\overline{g_{ns+t,2}}$ and $\lambda
_{ns+t,1}=\lambda_{ns+t,2}=\left\Vert f_{n}\right\Vert ^{-2}$. The integral
$\left\langle z^{t}f_{n}\left(  z^{s}\right)  ,z^{t}f_{n}\left(  z^{s}\right)
\right\rangle $ reduces to the case $s=1$ and $t=0$. When $s\geq1$ $\left\{
f_{n}^{0}\left(  z^{s}\right)  ,f_{n}^{1}\left(  z^{s}\right)  \right\}  $ is
an orthogonal basis for~$\mathcal{H}_{ns}$ and $\overline{z}^{s}f_{n-1}\left(
\overline{z}^{s}\right)  $ is orthogonal to $f_{n}\left(  z^{s}\right)  $. By
orthogonality $\left\Vert f_{n}\right\Vert ^{2}=\left\Vert f_{n}%
^{0}\right\Vert ^{2}+\left\Vert f_{n}^{1}\right\Vert ^{2}$ and the latter two
norms are standard Jacobi polynomial facts. The associated structural
constants are denoted by labeled $\lambda$'s. Thus
\begin{gather}
\lambda_{2n}^{0}    :=\left\Vert f_{2n}^{0}\right\Vert ^{-2}=\frac{n!\left(
\kappa_{0}+\kappa_{1}+1\right)  _{n}\left(  \kappa_{0}+\kappa_{1}+2n\right)
}{\left(  \kappa_{0}+\frac{1}{2}\right)  _{n}\left(  \kappa_{1}+\frac{1}
{2}\right)  _{n}\left(  \kappa_{0}+\kappa_{1}+n\right)  },\label{lb02n}\\
\lambda_{2n}^{1}    :=\left\Vert f_{2n}^{1}\right\Vert ^{-2}=\frac{\left(
n-1\right)  !\left(  \kappa_{0}+\kappa_{1}+1\right)  _{n}\left(  \kappa
_{0}+\kappa_{1}+2n\right)  }{\left(  \kappa_{0}+\frac{1}{2}\right)
_{n}\left(  \kappa_{1}+\frac{1}{2}\right)  _{n}},\label{lb12n}\\
\lambda_{2n}    :=\left\Vert f_{2n}\right\Vert ^{-2}=\frac{n!\left(
\kappa_{0}+\kappa_{1}+1\right)  _{n}}{\left(  \kappa_{0}+\frac{1}{2}\right)
_{n}\left(  \kappa_{1}+\frac{1}{2}\right)  _{n}}; \label{lb2n}
\end{gather}
and%
\begin{gather}
\lambda_{2n+1}^{0}    :=\left\Vert f_{2n+1}^{0}\right\Vert ^{-2}%
=\frac{n!\left(  \kappa_{0}+\kappa_{1}+1\right)  _{n}\left(  \kappa_{0}%
+\kappa_{1}+2n+1\right)  }{\left(  n+\kappa_{0}+\frac{1}{2}\right)  \left(
\kappa_{0}+\frac{1}{2}\right)  _{n+1}\left(  \kappa_{1}+\frac{1}{2}\right)
_{n+1}},\label{lb02n1}\\
\lambda_{2n+1}^{1}    :=\left\Vert f_{2n+1}^{1}\right\Vert ^{-2}%
=\frac{n!\left(  \kappa_{0}+\kappa_{1}+1\right)  _{n}\left(  \kappa_{0}%
+\kappa_{1}+2n+1\right)  }{\left(  n+\kappa_{1}+\frac{1}{2}\right)  \left(
\kappa_{0}+\frac{1}{2}\right)  _{n+1}\left(  \kappa_{1}+\frac{1}{2}\right)
_{n+1}},\label{lb12n1}\\
\lambda_{2n+1}    :=\left\Vert f_{2n+1}\right\Vert ^{-2}=\frac{n!\left(
\kappa_{0}+\kappa_{1}+1\right)  _{n}}{\left(  \kappa_{0}+\frac{1}{2}\right)
_{n+1}\left(  \kappa_{1}+\frac{1}{2}\right)  _{n+1}}. \label{lb2n1}%
\end{gather}
From this point on we no longer need the measure $\mu$ on the circle. Only the
algebraic expressions are used. The condition $\kappa_{0},\kappa_{1}\geq0$ is
replaced by the requirement that none of $-\kappa_{0}+\frac{1}{2}$, $-\kappa
_{1}+\frac{1}{2}$, $-s\left(  \kappa_{0}+\kappa_{1}\right)  $ equal a positive
integer. The exceptional case $-\left(  \kappa_{0}+\kappa_{1}\right)  \in%
\mathbb{N}
$ is taken up in the last section. In the next section we compute the
structural constants for the biorthogonal bases $\left\{  f_{n}\left(
z^{s}\right)  ,f_{n}\left(  \overline{z}^{s}\right)  \right\}  $ and $\left\{
z^{s}f_{n-1}\left(  z^{s}\right)  ,\overline{z}^{s}f_{n-1}\left(  \overline
{z}^{s}\right)  \right\}  $ for $\mathcal{H}_{ns}$ (see \cite[p.~461]{D2}. It
is easier to carry this out with material developed in the next section.

\section[Expressions for coefficients]{Expressions for coef\/f\/icients}

This is a detailed study of the coef\/f\/icients of $f_{n}\left(  z\right)  $ in
terms of powers of $z$, $\overline{z}$. The expressions are in the form of a
single sum of hypergeometric $_{3}F_{2}$-type, and can not be simplif\/ied any
further. For a polynomial $f$ in $z$, $\overline{z}$ def\/ine $c\left(
f;a,b\right)  $ to be the coef\/f\/icient of $z^{a}\overline{z}^{b}$ in $f$, that
is,
\begin{gather*}
f\left(  z,\overline{z}\right)  =\sum_{a,b\geq0}c\left(  f;a,b\right)
z^{a}\overline{z}^{b}.
\end{gather*}
Since we restrict to polynomials with real coef\/f\/icients the equation $c\big(
\overline{f\left(  z\right)  };a,b\big)  =c\left(  f;b,a\right)  $ is valid.
Further $c\left(  f\left(  z^{s}\right)  ;as,bs\right)  =c\left(
f;a,b\right)  $. Recall%
\[
K_{n}\left(  z,w\right)  :=\frac{1}{2^{n}n!}V^{z}\left(  \left(  z\overline
{w}+\overline{z}w\right)  ^{n}\right)  ,
\]
thus $Vz^{n-j}\overline{z}^{j}$ is $2^{n}j!\left(  n-j\right)  !$ times the
coef\/f\/icient of $w^{j}\overline{w}^{n-j}$ in $K_{n}\left(  z,w\right)  $. To
adapt the notation from equation (\ref{pkn}) for $P_{0}$ set $h_{01}%
=g_{01}=\lambda_{01}=1$ and $h_{02}=g_{02}=\lambda_{02}=0$. Then
\begin{gather}
K_{n}\left(  z,w\right)     =2^{-n}\sum_{j=0}^{\left\lfloor n/2\right\rfloor
}\frac{1}{j!\left(  s\kappa_{0}+s\kappa_{1}+1\right)  _{n-j}}\label{Kser}\\
\phantom{K_{n}\left(  z,w\right)     =}{} \times\left(  z\overline{z}w\overline{w}\right)  ^{j}\sum_{i=1}^{2}%
\lambda_{n-2j,i}h_{n-2j,i}\left(  z,\overline{z}\right)  g_{n-2j,i}\left(
\overline{w},w\right). \nonumber
\end{gather}

\begin{proposition}
For $0\leq m\leq n$,%
\begin{gather}
V\left(  z^{n-m}\overline{z}^{m}\right)     =m!\left(  n-m\right)
!\sum_{j=0}^{\left\lfloor n/2\right\rfloor }\frac{1}{j!\left(  s\kappa
_{0}+s\kappa_{1}+1\right)  _{n-j}}\label{vk}\\
\phantom{V\left(  z^{n-m}\overline{z}^{m}\right)     =}{}  \times\left(  z\overline{z}\right)  ^{j}\sum_{i=1}^{2}\lambda
_{n-2j,i}c\left(  g_{n-2j,i};n-m-j,m-j\right)  h_{n-2j,i}\left(
z,\overline{z}\right)  .\nonumber
\end{gather}
\end{proposition}

The nonzero terms appear at increments (in $j$) of $2s$. We start by f\/inding
$c\left(  f_{n}^{0};n-j,j\right)  $ and $c\left(  f_{n}^{1};n-j,j\right)  $.
This is straightforward and will serve as motivation for introducing a
specif\/ic useful $_{3}F_{2}$-series. Consider $f_{2n}^{0}\left(  z\right)  $
and recall that
\[
P_{n}^{\left(  \alpha,\beta\right)  }\left(  t\right)  =\frac{\left(
\alpha+1\right)  _{n}}{n!}~_{2}F_{1}\left(
\genfrac{}{}{0pt}{}{-n,n+\alpha+\beta+1}{\alpha+1}%
;\frac{1-t}{2}\right)  .
\]
When $z=re^{\mathrm{i}\theta}$ we have $\frac{1}{2}\left(  1-\cos
2\theta\right)  =-\left(  z-\overline{z}\right)  ^{2}/\left(  4z\overline
{z}\right)  $ so%
\begin{gather*}
r^{2n}P_{n}^{\left(  \kappa_{0}-\frac{1}{2},\kappa_{1}-\frac{1}{2}\right)
}\left(  \cos2\theta\right)  \\
\qquad{}=\frac{\left(  \kappa_{0}+\frac{1}{2}\right)  _{n}}{n!}\sum_{l=0}^{n}\sum
_{i=0}^{2l}\frac{\left(  -n\right)  _{l}\left(  n+\kappa_{0}+\kappa
_{1}\right)  _{l}\left(  2l\right)  !}{l!\left(  \kappa_{0}+\frac{1}%
{2}\right)  _{l}i!\left(  2l-i\right)  !}2^{-2l}z^{n+l-i}\overline{z}%
^{n-l+i}\left(  -1\right)  ^{l+i}\\
\qquad{}=\frac{\left(  \kappa_{0}+\frac{1}{2}\right)  _{n}}{n!}\sum_{j=-n}^{n}\left(
-1\right)  ^{j}z^{n+j}\overline{z}^{n-j}\sum_{i=\max\left(  -2j,0\right)
}^{n-j}\frac{\left(  -n\right)  _{j+i}\left(  n+\kappa_{0}+\kappa_{1}\right)
_{j+i}\left(  \frac{1}{2}\right)  _{j+i}}{i!\left(  \kappa_{0}+\frac{1}%
{2}\right)  _{j+i}\left(  2j+i\right)  !};
\end{gather*}
(substituting $l=i+j$, so $0\leq i+j\leq n$ and $0\leq i\leq2i+2j$ are the
ranges of the summation) by the $\left(  z,\overline{z}\right)  $-symmetry it
suf\/f\/ices to consider $j\geq0$. Thus%
\begin{gather*}
c\left(  f_{2n}^{0};n+j,n-j\right)  =\frac{\left(  \kappa_{0}+\frac{1}%
{2}\right)  _{n}\left(  -n\right)  _{j}\left(  n+\kappa_{0}+\kappa_{1}\right)
_{j}\left(  \frac{1}{2}\right)  _{j}}{\left(  \kappa_{0}+\frac{1}{2}\right)
_{j}\left(  2j\right)  !n!}\left(  -1\right)  ^{j}\\
\qquad{}\times\sum_{i=0}^{n-j}\frac{\left(  j-n\right)  _{i}\left(  n+\kappa
_{0}+\kappa_{1}+j\right)  _{i}\left(  \frac{1}{2}+j\right)  _{i}}{i!\left(
\kappa_{0}+\frac{1}{2}+j\right)  _{i}\left(  2j+1\right)  _{i}}\\
\qquad{}=\frac{\left(  n+\kappa_{0}+\kappa_{1}\right)  _{j}\left(  \kappa_{0}+\frac
{1}{2}+j\right)  _{n-j}}{2^{2j}j!\left(  n-j\right)  !}~_{3}F_{2}\left(
\genfrac{}{}{0pt}{}{j-n,n+\kappa_{0}+\kappa_{1}+j,j+\frac{1}{2}}{\kappa
_{0}+\frac{1}{2}+j,2j+1}%
;1\right)  ;
\end{gather*}
this used $\left(  2j\right)  !=2^{2j}j!\left(  \frac{1}{2}\right)  _{j}$ and
$\left(  \kappa_{0}+\frac{1}{2}\right)  _{n}/\left(  \kappa_{0}+\frac{1}%
{2}\right)  _{j}=\left(  \kappa_{0}+\frac{1}{2}+j\right)  _{n-j}$. The sum,
which appears to be a mysterious combination of the parameters, actually has a
nice form revealing more useful information.

\begin{definition}
For $n\in
\mathbb{N}
_{0}$ and parameters $a$, $b$, $c_{1}$, $c_{2}$ let
\begin{gather*}
E_{n}\left(  a,b;c_{1},c_{2}\right)     :=\frac{\left(  a\right)  _{n}\left(
c_{2}\right)  _{n}}{n!\left(  c_{1}+c_{2}\right)  _{n}}~_{3}F_{2}\left(
\genfrac{}{}{0pt}{}{-n,b,c_{1}}{1-n-a,1-c_{2}-n}%
;1\right) \\
\phantom{E_{n}\left(  a,b;c_{1},c_{2}\right) }{}  =\frac{1}{n!\left(  c_{1}+c_{2}\right)  _{n}}\sum_{j=0}^{n}\frac{\left(
-n\right)  _{j}}{j!}\left(  a\right)  _{n-j}\left(  b\right)  _{j}\left(
c_{1}\right)  _{j}\left(  c_{2}\right)  _{n-j}.
\end{gather*}

\end{definition}

Observe the symmetry $E_{n}\left(  a,b;c_{1},c_{2}\right)  =\left(  -1\right)
^{n}E_{n}\left(  b,a;c_{2},c_{1}\right)  $. This follows from manipulations
such as $\left(  a\right)  _{n-j}=\left(  -1\right)  ^{j}\left(  a\right)
_{n}/\left(  1-n-a\right)  _{j}$. The following transformation is relevant to
the calculation of coef\/f\/icients.

\begin{proposition}
\label{E3F2}For $n\in%
\mathbb{N}
_{0}$ and parameters $a,b,c_{1},c_{2}$%
\[
E_{n}\left(  a,b;c_{1},c_{2}\right)  =\frac{\left(  a+c_{1}\right)  _{n}}%
{n!}~_{3}F_{2}\left(
\genfrac{}{}{0pt}{}{-n,n+a+b+c_{1}+c_{2}-1,c_{1}}{a+c_{1},c_{1}+c_{2}}%
;1\right)  .
\]
\end{proposition}

\begin{proof}
Use the transformation%
\[
_{3}F_{2}\left(
\genfrac{}{}{0pt}{}{-n,A,B}{C,D}
;1\right)  =\frac{\left(  D-B\right)  _{n}}{\left(  D\right)  _{n}}~_{3}%
F_{2}\left(
\genfrac{}{}{0pt}{}{-n,C-A,B}{C,1+B-D-n}%
;1\right)  .
\]
First set $A=b,B=c_{1},C=1-n-a,D=1-n-c_{2}$ then%
\begin{gather*}
  \frac{\left(  a\right)  _{n}\left(  c_{2}\right)  _{n}}{n!\left(
c_{1}+c_{2}\right)  _{n}}~_{3}F_{2}\left(
\genfrac{}{}{0pt}{}{-n,b,c_{1}}{1-n-a,1-c_{2}-n}%
;1\right)   =\frac{\left(  a\right)  _{n}}{n!}~_{3}F_{2}\left(
\genfrac{}{}{0pt}{}{-n,1-n-a-b,c_{1}}{1-n-a,c_{1}+c_{2}}%
;1\right)  .
\end{gather*}
Set $A=1-n-a-b$, $B=c_{1}$, $C=c_{1}+c_{2}$, $D=1-n-a$ to obtain the stated formula. In
the calculation the reversal such as $\left(  1-n-a\right)  _{n}=\left(
-1\right)  ^{n}\left(  a\right)  _{n}$ is used several times.
\end{proof}

We arrive at a pleasing formula:
\[
c\left(  f_{2n}^{0};n+j,n-j\right)  =\frac{\left(  n+\kappa_{0}+\kappa
_{1}\right)  _{j}}{2^{2j}j!}E_{n-j}\left(  \kappa_{0},\kappa_{1};j+\frac{1}%
{2},j+\frac{1}{2}\right)  .
\]
It is useful because it clearly displays the result of setting one or both
parameters equal to zero (or a negative integer). That is $E_{n}\left(
\kappa_{0},0;c_{1},c_{2}\right)  =\frac{\left(  \kappa_{0}\right)  _{n}\left(
c_{2}\right)  _{n}}{n!\left(  c_{1}+c_{2}\right)  _{n}}$ and $n\geq1$ implies
$E_{n}\left(  0,0;c_{1},c_{2}\right)  =0$. (When $\kappa_{0}=\kappa_{1}=0$ the
polynomial $f_{2n}^{0}$ is a multiple of the Chebyshev polynomial of the f\/irst
kind, that is $f_{2n}^{0}\left(  z\right)  =\frac{\left(  n\right)  _{n}%
}{n!2^{2n}}\left(  z^{2n}+\overline{z}^{2n}\right)  $, a fact obvious from the
def\/inition of $f_{2n}^{0}$.) The remaining basis polynomials can all be
expressed in terms of the function $E$.

\begin{proposition}
\label{f2ncf}For $n\in \mathbb{N}$
\begin{gather*}
f_{2n}^{0}\left(  z\right)     =\sum_{j=1}^{n}\left(  z^{n+j}\overline
{z}^{n-j}+z^{n-j}\overline{z}^{n+j}\right)  \frac{1}{2^{2j}j!}\\
\phantom{f_{2n}^{0}\left(  z\right)     =}{}  \times\left(  n+\kappa_{0}+\kappa_{1}\right)  _{j}~E_{n-j}\left(
\kappa_{0},\kappa_{1};j+\frac{1}{2},j+\frac{1}{2}\right)
  +z^{n}\overline{z}^{n}E_{n}\left(  \kappa_{0},\kappa_{1};\frac{1}{2}%
,\frac{1}{2}\right)  ,\\
f_{2n}^{1}\left(  z\right)     =\sum_{j=1}^{n}\left(  z^{n+j}\overline
{z}^{n-j}-z^{n-j}\overline{z}^{n+j}\right)  \frac{1}{2^{2j}\left(  j-1\right)
!}\\
\phantom{f_{2n}^{1}\left(  z\right)     =}{}  \times\left(  n+\kappa_{0}+\kappa_{1}+1\right)  _{j-1}~E_{n-j}\left(
\kappa_{0},\kappa_{1};j+\frac{1}{2},j+\frac{1}{2}\right)  .
\end{gather*}

\end{proposition}

\begin{proof}
The expansion for $f_{2n}^{0}$ has already been determined. Next
$\frac{\mathrm{i}}{2}r^{2}\sin2\theta=\frac{1}{4}\left(  z^{2}-\overline
{z}^{2}\right)  $ so%
\begin{gather*}
\left(  \frac{\mathrm{i}}{2}r^{2}\sin2\theta\right)  r^{2n-2}P_{n-1}^{\left(
\kappa_{0}+\frac{1}{2},\kappa_{1}+\frac{1}{2}\right)  }\left(  \cos
2\theta\right)  =\frac{\left(  \kappa_{0}+\frac{3}{2}\right)  _{n-1}}{\left(
n-1\right)  !}\\
\qquad{}{}\times\sum_{l=0}^{n-1}\frac{\left(  1-n\right)  _{l}\left(  n+\kappa
_{0}+\kappa_{1}+1\right)  _{l}}{l!\left(  \kappa_{0}+\frac{3}{2}\right)  _{l}
}2^{-2l-2}\left(  -1\right)  ^{l}\left(  z+\overline{z}\right)  \left(
z-\overline{z}\right)  ^{2l+1}\left(  z\overline{z}\right)  ^{n-1-l},
\end{gather*}
and
\begin{gather*}
  2^{-2l-2}\left(  z+\overline{z}\right)  \left(  z-\overline{z}\right)
^{2l+1}\left(  z\overline{z}\right)  ^{n-1-l}\\
 \qquad{} =2^{-2l-2}\sum_{i=0}^{2l+2}\frac{\left(  2l+1\right)  !\left(
2l+2-2i\right)  }{i!\left(  2l+2-i\right)  !}\left(  -1\right)  ^{i}%
z^{n+l+1-i}\overline{z}^{n-l-1+i}\\
\qquad{}  =\sum_{i=0}^{2l+2}\frac{l!\left(  \frac{1}{2}\right)  _{l+1}\left(
l+1-i\right)  }{i!\left(  2l+2-i\right)  !}\left(  -1\right)  ^{i}%
z^{n+l+1-i}\overline{z}^{n-l-1+i};
\end{gather*}
substitute $l=j+i-1$. By the symmetry $\overline{f_{2n}^{1}\left(  z\right)
}=-f_{2n}^{1}\left(  z\right)  $ it suf\/f\/ices to f\/ind $c\left(  f_{2n}%
^{1};n+j,n-j\right)  $ for $1\leq j\leq n$. Indeed
\begin{gather*}
c\left(  f_{2n}^{1};n+j,n-j\right)     =\frac{j\left(  \kappa_{0}+\frac{3}%
{2}\right)  _{n-1}\left(  1-n\right)  _{j-1}\left(  n+\kappa_{0}+\kappa
_{1}+1\right)  _{j-1}\left(  \frac{1}{2}\right)  _{j}}{\left(  \kappa
_{0}+\frac{3}{2}\right)  _{j-1}\left(  2j\right)  !\left(  n-1\right)
!}\left(  -1\right)  ^{j-1}\\
\phantom{c\left(  f_{2n}^{1};n+j,n-j\right)     =}{}  \times\sum_{i=0}^{n-j}\frac{\left(  j-n\right)  _{i}\left(  n+\kappa
_{0}+\kappa_{1}+j\right)  _{i}\left(  \frac{1}{2}+j\right)  _{i}}{i!\left(
\kappa_{0}+\frac{1}{2}+j\right)  _{i}\left(  2j+1\right)  _{i}}\\
\phantom{c\left(  f_{2n}^{1};n+j,n-j\right)     }{}  =\frac{\left(  n+\kappa_{0}+\kappa_{1}+1\right)  _{j-1}}{2^{2j}\left(
j-1\right)  !}E_{n-j}\left(  \kappa_{0},\kappa_{1};j+\frac{1}{2},j+\frac{1}
{2}\right)  .
\end{gather*}
This completes the proof.
\end{proof}

\begin{proposition}
\label{f2n1cf}For $n\in
\mathbb{N}
_{0}$
\begin{gather*}
f_{2n+1}^{0}\left(  z\right)    =\left(  n+\kappa_{0}+\frac{1}{2}\right)
\sum_{j=0}^{n}\left(  z^{n+1+j}\overline{z}^{n-j}+z^{n-j}\overline{z}%
^{n+1+j}\right)  \frac{1}{2^{2j+1}j!}\\
\phantom{f_{2n+1}^{0}\left(  z\right)    =}{} \times\left(  n+\kappa_{0}+\kappa_{1}+1\right)  _{j}~E_{n-j}\left(
\kappa_{0},\kappa_{1};j+\frac{1}{2},j+\frac{3}{2}\right)  ,\\
f_{2n+1}^{1}\left(  z\right)     =\left(  n+\kappa_{1}+\frac{1}{2}\right)
\sum_{j=0}^{n}\left(  z^{n+1+j}\overline{z}^{n-j}-z^{n-j}\overline{z}%
^{n+1+j}\right)  \frac{1}{2^{2j+1}j!}\\
\phantom{f_{2n+1}^{1}\left(  z\right)     =}{}  \times\left(  n+\kappa_{0}+\kappa_{1}+1\right)  _{j}~E_{n-j}\left(
\kappa_{0},\kappa_{1};j+\frac{3}{2},j+\frac{1}{2}\right)  .
\end{gather*}
\end{proposition}

\begin{proof}
The second equation is straightforward:%
\begin{gather*}
f_{2n+1}^{1}\left(  z\right)  =\frac{1}{2}\left(  n+\kappa_{1}+\frac{1}%
{2}\right)  \left(  z-\overline{z}\right)  r^{2n}P_{n}^{\left(  \kappa
_{0}+\frac{1}{2},\kappa_{1}-\frac{1}{2}\right)  }\left(  \cos2\theta\right)
=\left(  n+\kappa_{1}+\frac{1}{2}\right)  \frac{\left(  \kappa_{0}+\frac{3}%
{2}\right)  _{n}}{n!}\\
\times\sum_{l=0}^{n}\sum_{i=0}^{2l+1}\frac{\left(  -n\right)  _{l}\left(
n+\kappa_{0}+\kappa_{1}+1\right)  _{l}\left(  2l+1\right)  !}{l!\left(
\kappa_{0}+\frac{3}{2}\right)  _{l}i!\left(  2l+1-i\right)  !}2^{-2l-1}%
z^{n+1+l-i}\overline{z}^{n-l+i}\left(  -1\right)  ^{l+i}\\
=\left(  n+\kappa_{1}+\frac{1}{2}\right)  \frac{\left(  \kappa_{0}+\frac{3}%
{2}\right)  _{n}}{n!}\sum_{j=0}^{n}\left(  z^{n+1+j}\overline{z}^{n-j}%
-z^{n-j}\overline{z}^{n+1-j}\right) \\
\times\left(  -1\right)  ^{j}\frac{\left(  -n\right)  _{j}\left(  n+\kappa
_{0}+\kappa_{1}+1\right)  _{j}\left(  \frac{1}{2}\right)  _{j+1}}{\left(
\kappa_{0}+\frac{3}{2}\right)  _{j}\left(  2j+1\right)  !}\sum_{i=0}%
^{n-j}\frac{\left(  j-n\right)  _{i}\left(  n+\kappa_{0}+\kappa_{1}%
+j+1\right)  _{i}\left(  \frac{3}{2}+j\right)  _{i}}{i!\left(  \kappa
_{0}+\frac{3}{2}+j\right)  _{i}\left(  2j+2\right)  _{i}},
\end{gather*}
(substituting $l=j+i$ for $0\leq j\leq n$) thus
\begin{gather*}
c\left(  f_{2n+1}^{1};n+1+j,n-j\right)  =-c\left(  f_{2n+1}^{1}%
;n-j,n+1+j\right) \\
\qquad{}=\left(  n+\kappa_{1}+\frac{1}{2}\right)  \frac{\left(  n+\kappa_{0}%
+\kappa_{1}+1\right)  _{j}}{j!2^{2j+1}}E_{n-j}\left(  \kappa_{0},\kappa
_{1};j+\frac{3}{2},j+\frac{1}{2}\right)  .
\end{gather*}
Note that $\left(  2j+1\right)  !=2^{2j+1}j!\left(  \frac{1}{2}\right)
_{j+1}$ and $\left(  \kappa_{0}+\frac{3}{2}\right)  _{n}/\left(  \kappa
_{0}+\frac{3}{2}\right)  _{j}=\left(  \kappa_{0}+\frac{3}{2}+j\right)  _{n-j}%
$. For $f_{2n+1}^{0}$ reverse the parameters, that is,
\[
P_{n}^{\left(  \kappa_{0}-\frac{1}{2},\kappa_{1}+\frac{1}{2}\right)  }\left(
\cos2\theta\right)  =\left(  -1\right)  ^{n}P_{n}^{\left(  \kappa_{1}+\frac
{1}{2},\kappa_{0}-\frac{1}{2}\right)  }\left(  -\cos2\theta\right)  ,
\]
and note $\frac{1}{2}\left(  1+\cos2\theta\right)  =\left(  z+\overline
{z}\right)  ^{2}/\left(  4z\overline{z}\right)  $ and $r\cos\theta=\frac{1}%
{2}\left(  z+\overline{z}\right)  $. Thus
\begin{gather*}
f_{2n+1}^{0}\left(  z\right)     =\left(  -1\right)  ^{n}\left(  n+\kappa
_{0}+\frac{1}{2}\right)  \frac{\left(  \kappa_{1}+\frac{3}{2}\right)  _{n}%
}{n!}\\
\phantom{f_{2n+1}^{0}\left(  z\right)}{}  \times\sum_{l=0}^{n}\sum_{i=0}^{2l+1}\frac{\left(  -n\right)  _{l}\left(
n+\kappa_{0}+\kappa_{1}+1\right)  _{l}\left(  2l+1\right)  !}{l!\left(
\kappa_{1}+\frac{3}{2}\right)  _{l}i!\left(  2l+1-i\right)  !}2^{-2l-1}%
z^{n+1+l-i}\overline{z}^{n-l+i}\\
\phantom{f_{2n+1}^{0}\left(  z\right)}{}  =\left(  -1\right)  ^{n}\left(  n+\kappa_{0}+\frac{1}{2}\right)
\frac{\left(  \kappa_{1}+\frac{3}{2}\right)  _{n}}{n!}\sum_{j=0}^{n}\left(
z^{n+1+j}\overline{z}^{n-j}+z^{n-j}\overline{z}^{n+1-j}\right) \\
\phantom{f_{2n+1}^{0}\left(  z\right)}{}  \times\frac{\left(  -n\right)  _{j}\left(  n+\kappa_{0}+\kappa
_{1}+1\right)  _{j}\left(  \frac{1}{2}\right)  _{j+1}}{\left(  \kappa
_{1}+\frac{3}{2}\right)  _{j}\left(  2j+1\right)  !}\sum_{i=0}^{n-j}%
\frac{\left(  j-n\right)  _{i}\left(  n+\kappa_{0}+\kappa_{1}+j+1\right)
_{i}\left(  \frac{3}{2}+j\right)  _{i}}{i!\left(  \kappa_{1}+\frac{3}%
{2}+j\right)  _{i}\left(  2j+2\right)  _{i}},
\end{gather*}
thus%
\begin{gather*}
c\left(  f_{2n+1}^{0};n+1+j,n-j\right)  =c\left(  f_{2n+1}^{0}%
;n-j,n+1+j\right) \\
\qquad{}=\left(  n+\kappa_{0}+\frac{1}{2}\right)  \frac{\left(  n+\kappa_{0}%
+\kappa_{1}+1\right)  _{j}}{j!2^{2j+1}}\left(  -1\right)  ^{n-j}E_{n-j}\left(
\kappa_{1},\kappa_{0};j+\frac{3}{2},j+\frac{1}{2}\right)  .
\end{gather*}
The symmetry relation $E_{m}\left(  b,a;c_{2},c_{1}\right)  =\left(
-1\right)  ^{m}E_{m}\left(  a,b;c_{1},c_{2}\right)  $ f\/inishes the computation.
\end{proof}

To f\/ind the coef\/f\/icients of $f_{n}$ we use contiguity relations satisf\/ied by
$E_{m}$.

\begin{lemma}
For $m\in
\mathbb{N}
_{0}$ and parameters $a$, $b$, $c$
\begin{gather}
  \left(  m+a+c\right)  E_{m}\left(  a,b;c,c+1\right)  -\left(  m+b+c\right)
E_{m}\left(  a,b;c+1,c\right) \label{eminus}\\
\qquad{} =2\left(  m+1\right)  E_{m+1}\left(  a,b;c,c\right)  ,\nonumber
\\
  \left(  m+a+c\right)  E_{m}\left(  a,b;c,c+1\right)  +\left(  m+b+c\right)
E_{m}\left(  a,b;c+1,c\right) \label{eplus}\\
\qquad{} =\frac{m+2c+1}{2c+1}\left(  m+a+b+2c\right)  E_{m}\left(
a,b;c+1,c+1\right)  .\nonumber
\end{gather}
\end{lemma}

\begin{proof}
We compute the coef\/f\/icient of $\left(  b\right)  _{j}$ for $0\leq j\leq m+1$
in the two identities. Note that $\left(  m+b+c\right)  \left(  b\right)
_{j}=\left(  b\right)  _{j+1}+\left(  m+c-j\right)  \left(  b\right)  _{j}$,
then replace $j$ by $j-1$ for the f\/irst term. The coef\/f\/icient of $\left(
b\right)  _{j}$ in $\left(  m+b+c\right)  E_{m}\left(  a,b;c+1,c\right)  $ is
\begin{gather*}
  \frac{1}{m!\left(  2c+1\right)  _{m}j!}\\
\times  \big\{  \left(  -m\right)  _{j}\left(  m+c-j\right)  \left(  a\right)
_{m-j}\left(  c\right)  _{m-j}\left(  c+1\right)  _{j}+j\left(  -m\right)
_{j-1}\left(  a\right)  _{m+1-j}\left(  c\right)  _{m+1-j}\left(  c+1\right)
_{j-1}\big\} \\
  =\frac{\left(  -m\right)  _{j-1}}{m!\left(  2c+1\right)  _{m}j!}\left(
a\right)  _{m-j}\left(  c\right)  _{m+1-j}\left(  c+1\right)  _{j-1}\left\{
\left(  -m+j-1\right)  \left(  c+j\right)  +j\left(  a+m-j\right)  \right\}  .
\end{gather*}
The coef\/f\/icient of $\left(  b\right)  _{j}$ in the left side of \eqref{eminus}
is
\begin{gather*}
  \frac{\left(  -m\right)  _{j-1}}{m!\left(  2c+1\right)  _{m}j!}\left(
a\right)  _{m-j}\left(  c\right)  _{j}\left(  c+1\right)  _{m-j}\\
 \qquad{} \times\left\{  \left(  m+a+c\right)  \left(  -m+j-1\right)  -\left(
-m+j-1\right)  \left(  c+j\right)  -j\left(  a+m-j\right)  \right\} \\
\qquad{}  =\frac{\left(  -m\right)  _{j-1}}{m!\left(  2c+1\right)  _{m}j!}\left(
a\right)  _{m-j}\left(  c\right)  _{j}\left(  c+1\right)  _{m-j}\left(
a+m-j\right)  \left(  -m-1\right) \\
\qquad{}  =\frac{2\left(  -1-m\right)  _{j}}{m!\left(  2c\right)  _{m+1}j!}\left(
a\right)  _{m+1-j}\left(  c\right)  _{j}\left(  c\right)  _{m+1-j}.
\end{gather*}
This proves equation \eqref{eminus}. For the right side of (\ref{eplus}) the
coef\/f\/icient of $\left(  b\right)  _{j}$ is found similarly as before ($\left(
m+a+b+2c\right)  \left(  b\right)  _{j}\allowbreak=\left(  m+a+2c-j\right)
\left(  b\right)  _{j}+\left(  b\right)  _{j+1}$, and so on). The coef\/f\/icient
of~$\left(  b\right)  _{j}$ in the left side is%
\[
\frac{\left(  -m\right)  _{j-1}}{m!\left(  2c+1\right)  _{m}j!}\left(
a\right)  _{m-j}\left(  c\right)  _{j}\left(  c+1\right)  _{m-j}\left\{
\left(  m+a+2c\right)  \left(  -m+j-1\right)  +j\left(  a-1\right)  \right\},
\]
and in the right side%
\begin{gather*}
  \frac{m+2c+1}{m!\left(  2c+1\right)  \left(  2c+2\right)  _{m}}%
\Bigg\{\frac{\left(  -m\right)  _{j}}{j!}\left(  m+a+2c-j\right)  \left(  a\right)
_{m-j}\left(  c+1\right)  _{j}\left(  c+1\right)  _{m-j}\\
\qquad{}  +\frac{\left(  -m\right)  _{j-1}}{j!}j\left(  a\right)  _{m+1-j}\left(
c+1\right)  _{j-1}\left(  c+1\right)  _{m+1-j}\Bigg\}\\
\qquad{}  =\frac{\left(  -m\right)  _{j-1}}{m!\left(  2c+1\right)  _{m}j!}\left(
a\right)  _{m-j}\left(  c+1\right)  _{j-1}\left(  c+1\right)  _{m-j}\\
\qquad{}  \times\left\{  \left(  -m+j-1\right)  \left(  m+a+2c-j\right)  \left(
c+j\right)  +j\left(  a+m-j\right)  \left(  c+m+1-j\right)  \right\}  ;
\end{gather*}
the expression in $\left\{  \cdot\right\}  $ equals $c\left(  m+a+2c\right)
\left(  -m+j-1\right)  +cj\left(  a-1\right)  $ which proves (\ref{eplus}).
\end{proof}

\begin{proposition}
\label{cofzfn}For $n\in
\mathbb{N}
_{0}$
\begin{gather*}
f_{2n}\left(  z\right)     =\sum_{j=1}^{n}\left(  \left(  n+\kappa_{0}%
+\kappa_{1}+j\right)  z^{n+j}\overline{z}^{n-j}+\left(  n+\kappa_{0}%
+\kappa_{1}-j\right)  z^{n-j}\overline{z}^{n+j}\right) \\
\phantom{f_{2n}\left(  z\right)     =}{}  \times\frac{1}{2^{2j}j!}\left(  n+\kappa_{0}+\kappa_{1}+1\right)
_{j-1}~E_{n-j}\left(  \kappa_{0},\kappa_{1};j+\frac{1}{2},j+\frac{1}{2}\right)
\\
\phantom{f_{2n}\left(  z\right)     =}{}  +E_{n}\left(  \kappa_{0},\kappa_{1};\frac{1}{2},\frac{1}{2}\right)
z^{n}\overline{z}^{n},
\\
f_{2n+1}\left(  z\right)     =\sum_{j=1}^{n+1}\left(  \left(  n+1+j\right)
z^{n+j}\overline{z}^{n+1-j}+\left(  n+1-j\right)  z^{n-j}\overline{z}%
^{n+1+j}\right) \\
\phantom{f_{2n+1}\left(  z\right)     =}{}  \times\frac{1}{2^{2j}j!}\left(  n+\kappa_{0}+\kappa_{1}+1\right)
_{j}~E_{n+1-j}\left(  \kappa_{0},\kappa_{1};j+\frac{1}{2},j+\frac{1}{2}\right)
\\
\phantom{f_{2n+1}\left(  z\right)     =}{}  +\left(  n+1\right)  E_{n+1}\left(  \kappa_{0},\kappa_{1};\frac{1}{2}%
,\frac{1}{2}\right)  z^{n}\overline{z}^{n+1}.
\end{gather*}
\end{proposition}

\begin{proof}
Recall $f_{n}=f_{n}^{0}+f_{n}^{1}$. For $0\leq j\leq n$ from Proposition
\ref{f2ncf} we f\/ind%
\begin{gather*}
c\left(  f_{2n};n+j,n-j\right)  =c\left(  f_{2n}^{0};n+j,n-j\right)  +c\left(
f_{2n}^{1};n+j,n-j\right) \\
\qquad{}=\frac{\left(  n+\kappa_{0}+\kappa_{1}+1\right)  _{j-1}}{2^{2j}j!}\left(
\left(  n+\kappa_{0}+\kappa_{1}\right)  +j\right)  E_{n-j}\left(  \kappa
_{0},\kappa_{1};j+\frac{1}{2},j+\frac{1}{2}\right)  ,\\
c\left(  f_{2n};n-j,n+j\right)  =c\left(  f_{2n}^{0};n+j,n-j\right)  -c\left(
f_{2n}^{1};n+j,n-j\right) \\
\qquad{}=\frac{\left(  n+\kappa_{0}+\kappa_{1}+1\right)  _{j-1}}{2^{2j}j!}\left(
\left(  n+\kappa_{0}+\kappa_{1}\right)  -j\right)  E_{n-j}\left(  \kappa
_{0},\kappa_{1};j+\frac{1}{2},j+\frac{1}{2}\right)  .
\end{gather*}
It remains to compute $c\left(  f_{2n+1};n+1+j,n-j\right)  $ and $c\left(
f_{2n+1};n-j,n+1-j\right)  $ for $0\leq j\leq n$. Write the arguments as
$\left(  n+\frac{1}{2}+\varepsilon\left(  j+\frac{1}{2}\right)  ,n+\frac{1}%
{2}-\varepsilon\left(  j+\frac{1}{2}\right)  \right)  $ with $\varepsilon
=\pm1$. Then, by Proposition \ref{f2n1cf},
\begin{gather*}
  c\left(  f_{2n+1};n+\frac{1}{2}+\varepsilon\left(  j+\frac{1}{2}\right)
,n+\frac{1}{2}-\varepsilon\left(  j+\frac{1}{2}\right)  \right) \\
\qquad{}  =c\left(  f_{2n+1}^{0};n+1+j,n-j\right)  +\varepsilon c\left(  f_{2n+1}%
^{1};n+1+j,n-j\right) \\
\qquad{}  =\frac{\left(  n+\kappa_{0}+\kappa_{1}+1\right)  _{j}}{j!2^{2j+1}}\Bigg\{\left(
n+\kappa_{0}+\frac{1}{2}\right)  E_{n-j}\left(  \kappa_{0},\kappa_{1}%
;j+\frac{1}{2},j+\frac{3}{2}\right) \\
\qquad{} +\varepsilon\left(  n+\kappa_{1}+\frac{1}{2}\right)  E_{n-j}\left(
\kappa_{0},\kappa_{1};j+\frac{3}{2},j+\frac{1}{2}\right)  \Bigg\}.
\end{gather*}
When $\varepsilon=1$ by \eqref{eplus} we obtain%
\begin{gather*}
c\left(  f_{2n+1};n+1+j,n-j\right)      =\frac{\left(  n+\kappa_{0}+\kappa_{1}+1\right)  _{j+1}}{\left(  j+1\right)
!2^{2j+2}}\left(  n+j+2\right)  E_{n-j}\left(  \kappa_{0},\kappa_{1}%
;j+\frac{3}{2},j+\frac{3}{2}\right)  ,
\end{gather*}
and when $\varepsilon=-1$ by \eqref{eminus} we obtain
\begin{gather*}
c\left(  f_{2n+1};n-j,n+1-j\right)  =
  \frac{\left(  n+\kappa_{0}+\kappa_{1}+1\right)  _{j}}{j!2^{2j}}\left(
n-j+1\right)  E_{n-j+1}\left(  \kappa_{0},\kappa_{1};j+\frac{1}{2},j+\frac
{1}{2}\right)  .
\end{gather*}
The stated formula for $f_{2n+1}$ uses $c\left(  f_{2n+1};n+j,n+1-j\right)  $
explicitly ($j$ is shifted by $1$).
\end{proof}

For $\mathcal{H}_{ns}$ with $n>0$ we intend to use both the orthogonal basis
$\left\{  f_{n}^{0}\left(  z^{s}\right)  ,f_{n}^{1}\left(  z^{s}\right)
\right\}  $ as well as the biorthogonal bases $\big\{  f_{n}\left(
z^{s}\right)  ,\overline{f_{n}\left(  z^{s}\right)  }\big\}  $ and $\big\{
z^{s}f_{n-1}\left(  z^{s}\right)  ,\overline{z}^{s}\overline{f_{n-1}\left(
z^{s}\right)  }\big\}  $. For the latter we need the value of $\nu
_{n}:=\left\langle f_{n}\left(  z^{s}\right)  ,z^{s}f_{n-1}\left(
z^{s}\right)  \right\rangle $. Instead of doing the integral directly we use
the two formulae for $P_{ns}\left(  z,w\right)  $, that is,%
\begin{gather*}
P_{ns}\left(  z,w\right)     =\lambda_{n}^{0}f_{n}^{0}\left(  z^{s}\right)
f_{n}^{0}\left(  w^{s}\right)  +\lambda_{n}^{1}f_{n}^{1}\left(  z^{s}\right)
\overline{f_{n}^{1}\left(  w^{s}\right)  }\\
\phantom{P_{ns}\left(  z,w\right)}{}  =\nu_{n}^{-1}\left(  f_{n}\left(  z^{s}\right)  \overline{w}^{s}%
f_{n-1}\left(  \overline{w}^{s}\right)  +f_{n}\left(  \overline{z}^{s}\right)
w^{s}f_{n-1}\left(  w^{s}\right)  \right)  .
\end{gather*}
From the coef\/f\/icients of $\overline{w}^{ns}$ in the equation we obtain%
\[
\lambda_{n}^{0}c\left(  f_{n}^{0};0,n\right)  f_{n}^{0}\left(  z^{s}\right)
+\lambda_{n}^{1}c\left(  f_{n}^{1};n,0\right)  f_{n}^{1}\left(  z^{s}\right)
=\nu_{n}^{-1}c\left(  f_{n-1};n-1,0\right)  f_{n}\left(  z^{s}\right)  .
\]
But $f_{n}=f_{n}^{0}+f_{n}^{1}$ so by the linear independence of $\left\{
f_{n}^{0},f_{n}^{1}\right\}  $ there are two equations for $c_{n}$ (one is
redundant). Thus%
\[
\nu_{n}=\frac{c\left(  f_{n-1};n-1,0\right)  }{\lambda_{n}^{0}c\left(
f_{n}^{0};0,n\right)  }=\frac{c\left(  f_{n-1};n-1,0\right)  }{\lambda_{n}%
^{1}c\left(  f_{n}^{1};n,0\right)  }.
\]
The calculation has two cases depending on $n$ being even or odd:%
\begin{gather}
\nu_{2n}    =\frac{2\left(  \kappa_{0}+\frac{1}{2}\right)  _{n}\left(
\kappa_{1}+\frac{1}{2}\right)  _{n}}{\left(  n-1\right)  !\left(  \kappa
_{0}+\kappa_{1}+1\right)  _{n-1}\left(  \kappa_{0}+\kappa_{1}+2n\right)
},\qquad n\geq1,\label{nu1}\\
\nu_{2n+1}    =\frac{2\left(  \kappa_{0}+\frac{1}{2}\right)  _{n+1}\left(
\kappa_{1}+\frac{1}{2}\right)  _{n+1}}{n!\left(  \kappa_{0}+\kappa
_{1}+1\right)  _{n}\left(  \kappa_{0}+\kappa_{1}+2n+1\right)  },\qquad n\geq0.
\label{nu2}%
\end{gather}

\section{The intertwining operator}

We describe $Vz^{a}\overline{z}^{b}$ for $a\geq b$. It is helpful to consider
the representations of $I_{2}\left(  2s\right)  $ since $V$ commutes with the
group action on polynomials. Since $z\overline{z}$ is invariant it suf\/f\/ices to
consider $\left(  z\overline{z}\right)  ^{b}z^{a-b},$ or $z^{m}$. The residue
of $m\operatorname{mod}2s$ is the determining factor. Suppose $m\equiv
j\operatorname{mod}2s$ and $j\neq0,s$. The representation of $I_{2}\left(
2s\right)  $ on $\mathrm{span}\left\{  z^{m},\overline{z}^{m}\right\}  $ is
irreducible and isomorphic to the one on $\mathrm{span}\left\{  z^{j}%
,\overline{z}^{j}\right\}  $ if $1\leq j<s$, and to the one on $\mathrm{span}%
\left\{  \overline{z}^{2s-j},z^{2s-j}\right\}  $ if $s<j<2s$. If
$m\equiv0\operatorname{mod}2s$ then $\mathrm{span}\left\{  z^{m},\overline
{z}^{m}\right\}  $ is the direct sum of the identity and determinant
representations (on $%
\mathbb{C}
1$ and $%
\mathbb{C}
\left(  z^{2s}-\overline{z}^{2s}\right)  $ respectively). If $m\equiv
s\operatorname{mod}2s$ then $\mathrm{span}\left\{  z^{m},\overline{z}%
^{m}\right\}  $ is the direct sum of the two representations realized on $%
\mathbb{C}
\left(  z^{s}-\overline{z}^{s}\right)  $ and $%
\mathbb{C}
\left(  z^{s}+\overline{z}^{s}\right)  $ (these are relative invariants).
Recall $P_{m}\left(  z,w\right)  =\sum\limits_{i=1}^{2}\lambda_{mi}h_{mi}\left(
z,\overline{z}\right)  g_{mi}\left(  \overline{w},w\right)  $ and equation
(\ref{vk}) shows that the nonzero terms in the expansion of $Vz^{a}%
\overline{z}^{b}$ occur only when the condition
\begin{equation}
c\left(  g_{a+b-2j,i};a-j,b-j\right)  \neq0 \label{nz}%
\end{equation}
is satisf\/ied. If $m\equiv0\operatorname{mod}s$ then $g_{mi}\left(
\overline{w},w\right)  $ is a polynomial in $w^{s},\overline{w}^{s}$ thus
(\ref{nz}) is equivalent to $a-j\equiv b-j\equiv0\operatorname{mod}s$, in
particular $a\equiv b\operatorname{mod}s$. In this case suppose $a=us+r\geq
b=vs+r$ with $0\leq r<s$. Set $b-j=\left(  v-k\right)  s$ then
$j=ks+r$, $a-j=\left(  u-k\right)  s$, $a+b-2j=\left(  u+v-2k\right)  s,$ and $0\leq
k\leq v\leq u$. We see that the nonzero terms occur for $P_{\left(
u+v-2k\right)  s}$ with $0\leq k\leq v$.

If $m\equiv t\operatorname{mod}s$ and $1\leq t<s$ then $g_{m1}\left(
w,\overline{w}\right)  =w^{t}f_{\left(  m-t\right)  /s}\left(  w^{s}%
,\overline{w}^{s}\right)  $ and (\ref{nz}) implies $a-j\equiv
t\operatorname{mod}s$, $b-j\equiv0\operatorname{mod}s$; further $g_{m2}\left(
w,\overline{w}\right)  =g_{m1}\left(  \overline{w},w\right)  $ and (\ref{nz})
implies $a-j\equiv0\operatorname{mod}s$, $b-j\equiv t\operatorname{mod}s$.

\begin{theorem}
\label{Vz1}Suppose $a-b\equiv t\operatorname{mod}s$, $1\leq t<s$ and $a>b$.
Let $b=vs+r$ with $v\geq0$ and $0\leq r<s$ and $a=us+r+t$, then%
\begin{gather*}
V\big(  z^{a}\overline{z}^{b}\big)     =a!b!\sum_{k=0}^{v}\frac{1}{\left(
ks+r\right)  !\left(  s\kappa_{0}+s\kappa_{1}+1\right)  _{a+\left(
v-k\right)  s}}\lambda_{u+v-2k}\\
\phantom{V\big(  z^{a}\overline{z}^{b}\big)     =}{}  \times c\left(  f_{u+v-2k};u-k,v-k\right)  \left(  z\overline{z}\right)
^{ks+r}z^{t}f_{u+v-2k}\left(  z^{s}\right) \\
\phantom{V\big(  z^{a}\overline{z}^{b}\big)     =}{}  +a!b!\sum_{k=1-\left\lfloor \left(  r+t\right)  /s\right\rfloor }^{v}%
\frac{1}{\left(  \left(  k-1\right)  s+r+t\right)  !\left(  s\kappa
_{0}+s\kappa_{1}+1\right)  _{b+\left(  u-k+1\right)  s}}\lambda_{u+v+1-2k}\\
\phantom{V\big(  z^{a}\overline{z}^{b}\big)     =}{}  \times c\left(  f_{u+v+1-2k};v-k,u-k+1\right)  \left(  z\overline
{z}\right)  ^{\left(  k-1\right)  s+r+t}\overline{z}^{s-t}f_{u+v+1-2k}\left(
\overline{z}^{s}\right)  .
\end{gather*}
\end{theorem}

\begin{proof}
Since $0<a-b=\left(  u-v\right)  s+t$ we have $u\geq v$. For the f\/irst part of
the series, corresponding to $i=1$ in $P_{ns+t}$ let $b-j=\left(  v-k\right)
s$ with $k\leq v$; then $j=b-\left(  v-k\right)  s=ks+r$, implying $k\geq0$.
Further $a-j=\left(  u-k\right)  s+t$, $a+b-2j=\left(  u+v-2k\right)
s+t=a+b-2r-2ks$ (and $a+b-j=a+\left(  b-j\right)  =a+\left(  v-k\right)  s$).
Also $c\left(  z^{t}f_{u+v-2k}\left(  z^{s}\right)  ;a-j,b-j\right)  =c\left(
f_{u+v-2k};u-k,v-k\right)  .$ This proves the f\/irst part. For the second part,
with $i=2$ in $P_{ns-t}$ let $b-j=\left(  v-k\right)  s+\left(  s-t\right)  $,
thus $k\leq v$. Then $j=\left(  k-1\right)  s+r+t$. The requirement $j\geq0$
implies $1-k\leq\frac{r+t}{s}$, that is $k\geq1-\left\lfloor \frac{r+t}%
{s}\right\rfloor $ (if $0\leq r+t<s$ then $k\geq1$, otherwise $s\leq r+t<2s$
and $k\geq0$). Also $a-j=\left(  u-k+1\right)  s$ and $a+b-2j=\left(
u+v+1-2k\right)  s+\left(  s-t\right)  $ (and $a+b-j=b+\left(  a-j\right)
=b+\left(  u-k+1\right)  s$). In this case we use $c\left(  \overline{z}%
^{s-t}f_{u+v+1-2k}\left(  \overline{z}^{s}\right)  ;a-j,b-j\right)  =c\left(
f_{u+v+1-2k};v-k,u-k+1\right)  $.
\end{proof}

Note that the degrees of $f_{m}$ have the same parity as $u+v$ in the f\/irst
sum, and the opposite in the second sum. By Proposition~\ref{cofzfn} we can
f\/ind the coef\/f\/icients explicitly. If $u+v$ is even then
\begin{gather*}
c\left(  f_{u+v-2k};u-k,v-k\right)     =\frac{1}{2^{u-v}\left(  \frac{u-v}%
{2}\right)  !}\left(  \frac{u+v}{2}-k+\kappa_{0}+\kappa_{1}\right)
_{\frac{u-v}{2}}\\
\phantom{c\left(  f_{u+v-2k};u-k,v-k\right)     =}{}  \times E_{v-k}\left(  \kappa_{0},\kappa_{1};\frac{u-v+1}{2},\frac{u-v+1}%
{2}\right)  ,
\end{gather*}
and%
\begin{gather*}
c\left(  f_{u+v+1-2k};v-k,u-k+1\right)     =\frac{\left(  v-k+1\right)
}{2^{u-v}\left(  \frac{u-v}{2}\right)  !}\left(  \frac{u+v}{2}-k+\kappa
_{0}+\kappa_{1}+1\right)  _{\frac{u-v}{2}}\\
\phantom{c\left(  f_{u+v+1-2k};v-k,u-k+1\right)  =}{}  \times E_{v-k+1}\left(  \kappa_{0},\kappa_{1};\frac{u-v+1}{2},\frac
{u-v+1}{2}\right)  .
\end{gather*}
If $u+v$ is odd then%
\begin{gather*}
c\left(  f_{u+v-2k};u-k,v-k\right)     =\frac{\left(  u-k+1\right)
}{2^{u-v+1}\left(  \frac{u-v+1}{2}\right)  !}\left(  \frac{u+v+1}{2}%
-k+\kappa_{0}+\kappa_{1}\right)  _{\frac{u-v+1}{2}}\\
\phantom{c\left(  f_{u+v-2k};u-k,v-k\right)     =}{}  \times E_{v-k}\left(  \kappa_{0},\kappa_{1};\frac{u-v}{2}+1,\frac{u-v}%
{2}+1\right)  ,
\end{gather*}
and%
\begin{gather*}
c\left(  f_{u+v+1-2k};v-k,u-k+1\right)     =\frac{\left(  v-k+1\right)
}{2^{u-v+1}\left(  \frac{u-v+1}{2}\right)  !}\left(  \frac{u+v+3}{2}%
-k+\kappa_{0}+\kappa_{1}\right)  _{\frac{u-v-1}{2}}\\
\phantom{c\left(  f_{u+v+1-2k};v-k,u-k+1\right)     =}{}  \times E_{v-k}\left(  \kappa_{0},\kappa_{1};\frac{u-v}{2}+1,\frac{u-v}%
{2}+1\right)  .
\end{gather*}

\begin{theorem}
\label{Vz2}Suppose $a\equiv b\operatorname{mod}s$, and $a\geq b$. Let
$a=us+r\geq b=vs+r$ with $0\leq r<s$ and $v\geq0$. If $a>b$ then
\begin{gather*}
V\big(  z^{a}\overline{z}^{b}\big)  =a!b!\sum_{k=0}^{v}\frac{1}{\left(
ks+r\right)  !\left(  s\kappa_{0}+s\kappa_{1}+1\right)  _{b+\left(
u-k\right)  s}}\nu_{u+v-2k}^{-1}\\
\phantom{V\big(  z^{a}\overline{z}^{b}\big)  =}{}\times\left(  z\overline{z}\right)  ^{ks+r}\{c\left(  f_{u+v-2k-1}%
;u-k-1,v-k\right)  f_{u+v-2k}\left(  z^{s}\right) \\
\phantom{V\big(  z^{a}\overline{z}^{b}\big)  =}{}+c\left(  f_{u+v-2k-1};v-k-1,u-k\right)  f_{u+v-2k}\left(  \overline{z}%
^{s}\right)  \},\\
V\left(  \frac{1}{2}\big(  z^{a}\overline{z}^{b}-z^{b}\overline{z}%
^{a}\big)  \right)  =a!b!\sum_{k=0}^{v}\frac{1}{\left(  ks+r\right)
!\left(  s\kappa_{0}+s\kappa_{1}+1\right)  _{b+\left(  u-k\right)  s}}%
\lambda_{u+v-2k}^{1}\\
\phantom{V\left(  \frac{1}{2}\big(  z^{a}\overline{z}^{b}-z^{b}\overline{z}%
^{a}\big)  \right)  =}{}\times c\left(  f_{u+v-2k}^{1};u-k,v-k\right)  \left(  z\overline{z}\right)
^{ks+r}f_{u+v-2k}^{1}\left(  z^{s}\right)  .
\end{gather*}
If $a\geq b$ then
\begin{gather*}
V\left(  \frac{1}{2}\big(  z^{a}\overline{z}^{b}+z^{b}\overline{z}%
^{a}\big)  \right)     =a!b!\sum_{k=0}^{v}\frac{1}{\left(  ks+r\right)
!\left(  s\kappa_{0}+s\kappa_{1}+1\right)  _{b+\left(  u-k\right)  s}}%
\lambda_{u+v-2k}^{0}\\
\phantom{V\left(  \frac{1}{2}\big(  z^{a}\overline{z}^{b}+z^{b}\overline{z}%
^{a}\big)  \right)     =}{}  \times c\left(  f_{u+v-2k}^{0};u-k,v-k\right)  \left(  z\overline
{z}\right)  ^{ks+r}f_{u+v-2k}^{0}\left(  z^{s}\right)  .
\end{gather*}
\end{theorem}

\begin{proof}
The three dif\/ferent expansions for $z^{a}\overline{z}^{b}$, $\frac{1}{2}\left(
z^{a}\overline{z}^{b}-z^{b}\overline{z}^{a}\right)  $ and $\frac{1}{2}\left(
z^{a}\overline{z}^{b}+z^{b}\overline{z}^{a}\right)  $ use the bases $\left\{
f_{j},\overline{f_{j}}\right\}$, $\big\{  f_{j}^{1}\big\}  $ and $\big\{
f_{j}^{0}\big\}  $ respectively. Suppose $a=us+r\geq b=vs+r$ with $0\leq
r<s$. Set $b-j=\left(  v-k\right)  s$ then $j=ks+r$, $a-j=\left(  u-k\right)
s$, $a+b-2j=\left(  u+v-2k\right)  s,$ and $0\leq k\leq v\leq u$. Consider the
case $a>b$, that is, $u>v$. For arbitrary $m\geq1$ the basis $\left\{
f_{m}\left(  z^{s},\overline{z}^{s}\right)  ,f_{m}\left(  \overline{z}%
^{s},z^{s}\right)  \right\}  $ for $\mathcal{H}_{sm}$ has the biorthogonal set
$\left\{  z^{s}f_{m-1}\left(  z^{s},\overline{z}^{s}\right)  ,\overline{z}%
^{s}f_{m-1}\left(  \overline{z}^{s},z^{s}\right)  \right\}  $ and%
\begin{gather*}
c\left(  z^{s}f_{n_{1}+n_{2}-1}\left(  z^{s},\overline{z}^{s}\right)
;n_{1}s,n_{2}s\right)     =c\left(  f_{n_{1}+n_{2}-1};n_{1}-1,n_{2}\right)
,\\
c\left(  \overline{z}^{s}f_{n_{1}+n_{2}-1}\left(  \overline{z}^{s}%
,z^{s}\right)  ;n_{1}s,n_{2}s\right)     =c\left(  f_{n_{1}+n_{2}-1}%
;n_{2}-1,n_{1}\right)  .
\end{gather*}
The constants $\nu_{m}$ are given in equations (\ref{nu1}) and (\ref{nu2}).
This demonstrates the f\/irst series. The remaining two follow from Proposition~\ref{vk}.
\end{proof}

Observe that in the series for $V\left(  z^{a}\overline{z}^{b}\right)  $ the
lowest-degree term with $k=v<u$ reduces to one summand since $c\left(
f_{u-v-1};-1,u-v\right)  =0$. Each term in $V\left(  z^{a}\overline{z}%
^{b}-z^{b}\overline{z}^{a}\right)  $ is of the same representation type, $%
\mathbb{C}
\left(  z^{s}-\overline{z}^{s}\right)  $ when $a-b\equiv s\operatorname{mod}%
2s$ or $%
\mathbb{C}
\left(  z^{2s}-\overline{z}^{2s}\right)  $ when $a-b\equiv0\operatorname{mod}%
2s$. Similarly each term in $V\left(  z^{a}\overline{z}^{b}+z^{b}\overline
{z}^{a}\right)  $ is of the representation type $%
\mathbb{C}
\left(  z^{s}+\overline{z}^{s}\right)  $ or $%
\mathbb{C}
1$ (depending on the parity of $\frac{a-b}{s}$). The coef\/f\/icients can be found
from Propositions \ref{f2ncf} and \ref{f2n1cf}.

For $a>b$ consider $z^{a}\overline{z}^{b}$ as $\left(  z\overline{z}\right)
^{b}$ times the (ordinary) harmonic polynomial $z^{a-b}$. The fact that
$V\left(  z^{a}\overline{z}^{b}\right)  $ is $L^{2}\left(  \mathbb{T}%
,\mu\right)  $-orthogonal to $\mathcal{H}_{n}$ for $n<a-b$, equivalently, that
the above series for $V\left(  z^{a}\overline{z}^{b}\right)  $ contain no
terms involving $\mathcal{H}_{n}$ with $n<a-b$ (that is, a term like
$c_{n}\left(  z\overline{z}\right)  ^{\left(  a+b-n\right)  /2}p_{n}\left(
z\right)  $ with $p_{n}\in\mathcal{H}_{n}$), is a special case of a result of
Xu \cite{X}. This paper also has formulae for $Vz^{2m}$ when $s=2$, that is,
the group $I_{2}\left(  4\right)  $.

\section{Singular values}

The term \textquotedblleft singular values\textquotedblright\ refers to the
set $K^{s}$ of pairs $\left(  \kappa_{0},\kappa_{1}\right)  \in%
\mathbb{C}
^{2}$ for which $V$ is not def\/ined on all polynomials in $z,\overline{z}$. Let%
\[
K_{0}:=\left\{  \left(  \kappa_{0},\kappa_{1}\right)  \in%
\mathbb{C}
^{2}:\left\{  \kappa_{0},\kappa_{1}\right\}  \cap\left(  -\tfrac{1}{2}-%
\mathbb{N}
_{0}\right)  \neq\varnothing\right\}  ,
\]
(at least one of $\kappa_{0},\kappa_{1}$ is in $\left\{  -\frac{1}{2}%
,-\frac{3}{2},\ldots\right\}  $). It was shown by de Jeu, Opdam and the author
\cite[p.~248]{DJO} that $K^{s}=K_{0}\cup\left\{  \left(  \kappa_{0},\kappa
_{1}\right)  :\kappa_{0}+\kappa_{1}=-\frac{j}{s},j\in%
\mathbb{N}
,\frac{j}{s}\notin%
\mathbb{N}
\right\}  $. To illustrate how the singular values appear in the formulae for
$V$ consider $Vz^{2ns+1}$ (for $s>1,n\geq1$) which has only one term in the
formula from Theorem \ref{Vz1}. In particular
\[
c\left(  Vz^{2ns+1};2ns+1,0\right)  =\frac{\left(  2ns+1\right)  !\left(
\kappa_{0}+\kappa_{1}+1\right)  _{2n}\left(  n+\kappa_{0}+\kappa_{1}+1\right)
_{n}}{2^{4n}n!\left(  \kappa_{0}+\frac{1}{2}\right)  _{n}\left(  \kappa
_{1}+\frac{1}{2}\right)  _{n}\left(  s\kappa_{0}+s\kappa_{1}+1\right)
_{2ns+1}}.
\]
The denominator vanishes for $\kappa_{0},\kappa_{1}=-\frac{1}{2},-\frac{3}%
{2},\ldots,-\frac{2n-1}{2}$ and $\kappa_{0}+\kappa_{1}=-\frac{k}{s}$ for
$1\leq k\leq2ns+1$. There appear to be singularities at $\kappa_{0}+\kappa
_{1}=-k$ for $1\leq k\leq2n$ but the term $\left(  \kappa_{0}+\kappa
_{1}+1\right)  _{2n}$ in the numerator cancels these zeros. The same
cancellation occurs for arbitrary $V\left(  z^{a}\overline{z}^{b}\right)  $ in
a more complicated way. The formula for $K_{n}\left(  z,w\right)  $ has the
factors $\left(  s\left(  \kappa_{0}+\kappa_{1}\right)  +1\right)  _{n-j}$ in
the denominators thus the individual terms can have simple poles at
$\kappa_{0}+\kappa_{1}=-\frac{k}{s}$ for $k\in%
\mathbb{N}
$. We will show directly that the singularities at $\kappa_{0}+\kappa_{1}=-m$
are removable when $K_{n}\left(  z,w\right)  $ is expressed as a quotient of
polynomials in $\kappa_{0}$, $\kappa_{1}$. It turns out that the terms with poles
can be paired in such a way that the sum of each pair has a removable
singularity. The pairs correspond to $\left\{  P_{k},P_{2sm-k}\right\}  $ for
certain values of $k$.

Throughout we assume that $\left(  \kappa_{0},\kappa_{1}\right)  \notin K_{0}$.

We use an elementary algebraic result: suppose a rational function $F\left(
\alpha,\beta\right)  $ (with coef\/f\/icients in the ring $%
\mathbb{Q}
\left[  z,\overline{z},w,\overline{w}\right]  $) vanishes for a countable set
of values $\left\{  \alpha=0,\beta=r_{n}:n\in%
\mathbb{N}
_{0}\right\}  $ (which are not poles) then $F\left(  \alpha,\beta\right)  $ is
divisible by $\alpha$; indeed the numerator of $F\left(  0,\beta\right)  $ is
a polynomial in $\beta$ vanishing at all $\beta=r_{n}$ hence is zero. This
result will be applied with $\alpha=\kappa_{0}+\kappa_{1}+m$, $\beta=\kappa
_{0}-\kappa_{1}$.

Most of the section concerns the proof of the following result: let
$\kappa_{0}+\kappa_{1}=-m$ then $P_{N}\left(  z,w\right)  =0$ for $N>2sm$ and
$P_{N}\left(  z,w\right)  +\left(  z\overline{z}w\overline{w}\right)
^{N-sm}P_{2sm-N}\left(  z,w\right)  =0$ for $0\leq N\leq2sm$. The Poisson
kernels $P_{n}$ were described in equation (\ref{pkn}). There are a number of
cases, roughly corresponding to the representations of $I_{2}\left(
2s\right)  $.

\begin{proposition}
\label{f2n0m}Suppose $-\left(  \kappa_{0}+\kappa_{1}\right)  =m\in%
\mathbb{N}
$ then%
\begin{gather*}
f_{2n}^{0}\left(  z\right)     =\frac{\left(  \kappa_{0}+\frac{1}{2}\right)
_{n}\left(  m-n\right)  !}{\left(  \kappa_{0}+\frac{1}{2}\right)  _{m-n}%
n!}\left(  z\overline{z}\right)  ^{2n-m}f_{2m-2n}^{0}\left(  z\right)  ,\qquad 0\leq
n\leq m,\\
f_{2n}^{1}\left(  z\right)     =\frac{\left(  \kappa_{0}+\frac{1}{2}\right)
_{n}\left(  m-n-1\right)  !}{\left(  \kappa_{0}+\frac{1}{2}\right)
_{m-n}\left(  n-1\right)  !}\left(  z\overline{z}\right)  ^{2n-m}f_{2m-2n}%
^{1}\left(  z\right)  ,\qquad 1\leq n\leq m-1.
\end{gather*}
\end{proposition}

\begin{proof}
The argument uses the Jacobi polynomials directly. Recall $z{=}re^{\mathrm{i}%
\theta}$. Then for $0\leq n\leq m$
\begin{gather*}
f_{2n}^{0}\left(  z\right)     =r^{2n}\frac{\left(  \kappa_{0}+\frac{1}%
{2}\right)  _{n}}{n!}~_{2}F_{1}\left(
\genfrac{}{}{0pt}{}{-n,n-m}{\kappa_{0}+\frac{1}{2}}%
;\frac{1-\cos2\theta}{2}\right)  ,\\
f_{2m-2n}^{0}\left(  z\right)     =r^{2m-2n}\frac{\left(  \kappa_{0}+\frac
{1}{2}\right)  _{m-n}}{\left(  m-n\right)  !}~_{2}F_{1}\left(
\genfrac{}{}{0pt}{}{-\left(  m-n\right)  ,\left(  m-n\right)  -m}{\kappa
_{0}+\frac{1}{2}}%
;\frac{1-\cos2\theta}{2}\right)  ,
\end{gather*}
while for $1\leq n\leq m-1$
\begin{gather*}
f_{2n}^{1}\left(  z\right)     =\mathrm{i}r^{2n}\sin2\theta\frac{\left(
\kappa_{0}+\frac{3}{2}\right)  _{n-1}}{\left(  n-1\right)  !}~_{2}F_{1}\left(
\genfrac{}{}{0pt}{}{1-n,n-m+1}{\kappa_{0}+\frac{3}{2}}%
;\frac{1-\cos2\theta}{2}\right)  ,\\
f_{2m-2n}^{1}\left(  z\right)     =\mathrm{i}r^{2m-2n}\sin2\theta
\frac{\left(  \kappa_{0}+\frac{3}{2}\right)  _{m-n-1}}{\left(  m-n-1\right)
!}~_{2}F_{1}\left(
\genfrac{}{}{0pt}{}{-\left(  m-n-1\right)  ,1-n}{\kappa_{0}+\frac{1}{2}}%
;\frac{1-\cos2\theta}{2}\right)  .
\end{gather*}
This proves the formulae.
\end{proof}

\begin{proposition}
\label{f2n1m}Suppose $-\left(  \kappa_{0}+\kappa_{1}\right)  =m\in%
\mathbb{N}
$ and $0\leq n<m$ then%
\begin{gather*}
f_{2n+1}^{0}\left(  z\right)     =\frac{\left(  \kappa_{0}+\frac{1}%
{2}\right)  _{n+1}\left(  m-n-1\right)  !}{\left(  \kappa_{0}+\frac{1}%
{2}\right)  _{m-n}n!}\left(  z\overline{z}\right)  ^{2n-m+1}f_{2m-2n-1}%
^{0}\left(  z\right)  ,\\
f_{2n+1}^{1}\left(  z\right)     =\frac{\left(  \kappa_{0}+\frac{1}%
{2}\right)  _{n}\left(  m-n-1\right)  !}{\left(  \kappa_{0}+\frac{1}%
{2}\right)  _{m-n-1}n!}\left(  z\overline{z}\right)  ^{2n-m+1}f_{2m-2n-1}%
^{1}\left(  z\right)  .
\end{gather*}
\end{proposition}

\begin{proof}
Similarly to the even case we have%
\begin{gather*}
f_{2n+1}^{0}\left(  z\right)     =r^{2n+1}\cos\theta\frac{\left(  \kappa
_{0}+\frac{1}{2}\right)  _{n+1}}{n!}~_{2}F_{1}\left(
\genfrac{}{}{0pt}{}{-n,n-m+1}{\kappa_{0}+\frac{1}{2}}%
;\frac{1-\cos2\theta}{2}\right)  ,\\
f_{2m-2n-1}^{0}\left(  z\right)     =r^{2m-2n-1}\cos\theta\frac{\left(
\kappa_{0}+\frac{1}{2}\right)  _{m-n}}{\left(  m-n-1\right)  !}~_{2}%
F_{1}\left(
\genfrac{}{}{0pt}{}{-\left(  m-n-1\right)  ,-n}{\kappa_{0}+\frac{1}{2}}%
;\frac{1-\cos2\theta}{2}\right)  ,
\end{gather*}
and%
\begin{gather*}
f_{2n+1}^{1}\left(  z\right)     =\mathrm{i}r^{2n+1}\sin\theta\frac{\left(
\kappa_{1}+n+\frac{1}{2}\right)  \left(  \kappa_{0}+\frac{3}{2}\right)  _{n}%
}{n!}~_{2}F_{1}\left(
\genfrac{}{}{0pt}{}{-n,n-m+1}{\kappa_{0}+\frac{3}{2}}%
;\frac{1-\cos2\theta}{2}\right)  ,\\
f_{2m-2n-1}^{1}\left(  z\right)     =\mathrm{i}r^{2m-2n-1}\sin\theta
\frac{\left(  \kappa_{1}+m-n-\frac{1}{2}\right)  \left(  \kappa_{0}+\frac
{3}{2}\right)  _{m-n-1}}{\left(  m-n-1\right)  !}\\
\phantom{f_{2m-2n-1}^{1}\left(  z\right)     =}{} \times~_{2}F_{1}\left(
\genfrac{}{}{0pt}{}{-\left(  m-n-1\right)  ,-n}{\kappa_{0}+\frac{3}{2}}%
;\frac{1-\cos2\theta}{2}\right)  .
\end{gather*}
Thus
\begin{gather*}
\dfrac{f_{2n+1}^{1}\left(  z\right)  }{f_{2m-2n-1}^{1}\left(  z\right)
}
 =  r^{4n-2m+2}\dfrac{\left(  m-n-1\right)  !\left(  \kappa_{0}+\frac{1}%
{2}\right)  _{n+1}\left(  -m-\kappa_{0}+n+\frac{1}{2}\right)  }{n!\left(
\kappa_{0}+\frac{1}{2}\right)  _{m-n}\left(  -\kappa_{0}-n-\frac{1}{2}\right)
}
\\
\phantom{\dfrac{f_{2n+1}^{1}\left(  z\right)  }{f_{2m-2n-1}^{1}\left(  z\right)
}}{} =  r^{4n-2m+2}\dfrac{\left(  m-n-1\right)  !\left(  \kappa_{0}%
+\frac{1}{2}\right)  _{n}}{n!\left(  \kappa_{0}+\frac{1}{2}\right)  _{m-n-1}}.\tag*{\qed}
\end{gather*}\renewcommand{\qed}{}
\end{proof}

\begin{proposition}
\label{f2nm1}Suppose $-\left(  \kappa_{0}+\kappa_{1}\right)  =m\in%
\mathbb{N}
$ and $0\leq n<m$ then%
\[
f_{2n}\left(  z\right)  =\frac{\left(  \kappa_{0}+\frac{1}{2}\right)
_{n}\left(  m-n-1\right)  !}{\left(  \kappa_{0}+\frac{1}{2}\right)  _{m-n}%
n!}\left(  z\overline{z}\right)  ^{2n-m}\overline{z}f_{2m-2n-1}\left(
\overline{z}\right)  .
\]

\end{proposition}

\begin{proof}
We use the expressions from Proposition \ref{cofzfn}. First we show for $0\leq
j\leq\min\left(  n,m-n\right)  $ that
\[
E_{n-j}\left(  \kappa_{0},\kappa_{1};j+\frac{1}{2},j+\frac{1}{2}\right)
=\frac{\left(  \kappa_{0}+\frac{1}{2}\right)  _{m-n}\left(  n-j\right)
!}{\left(  \kappa_{0}+\frac{1}{2}\right)  _{n}\left(  m-n-j\right)
!}E_{m-n-j}\left(  \kappa_{0},\kappa_{1};j+\frac{1}{2},j+\frac{1}{2}\right)
.
\]
Indeed by Proposition \ref{E3F2}%
\begin{gather*}
  \left(  \kappa_{0}+\frac{1}{2}\right)  _{j}E_{n-j}\left(  \kappa_{0}%
,\kappa_{1};j+\frac{1}{2},j+\frac{1}{2}\right)
  =\frac{\left(  \kappa_{0}+\frac{1}{2}\right)  _{n}}{\left(  n-j\right)
!}~_{3}F_{2}\left(
\genfrac{}{}{0pt}{}{j-n,n-m+j,j+\frac{1}{2}}{\kappa_{0}+j+\frac{1}{2},2j+1}%
;1\right)  ,
\end{gather*}
and%
\begin{gather*}
  \left(\!  \kappa_{0}+\frac{1}{2}\!\right)  _{j}E_{m-n-j}\left(\!  \kappa
_{0},\kappa_{1};j+\frac{1}{2},j+\frac{1}{2}\!\right)
  =\frac{\left(  \kappa_{0}+\frac{1}{2}\right)  _{m-n}}{\left(  m-n-j\right)
!}~_{3}F_{2}\!\left(\!
\genfrac{}{}{0pt}{}{n-m+j,j-n,j+\frac{1}{2}}{\kappa_{0}+j+\frac{1}{2},2j+1}%
;1\!\right) \! .\!
\end{gather*}
Let
\begin{gather*}
g_{1}\left(  z\right)     :=\frac{n!}{\left(  \kappa_{0}+\frac{1}{2}\right)
_{n}}f_{2n}\left(  z\right)  ,\qquad
g_{2}\left(  z\right)     :=\frac{\left(  m-n-1\right)  !}{\left(  \kappa
_{0}+\frac{1}{2}\right)  _{m-n}}\left(  z\overline{z}\right)  ^{2n-m}%
\overline{z}f_{2m-2n-1}\left(  \overline{z}\right)  ,
\end{gather*}
and for $j\geq0$ let%
\[
b_{j}:=\frac{1}{2^{2j}j!}~_{3}F_{2}\left(
\genfrac{}{}{0pt}{}{n-m+j,j-n,j+\frac{1}{2}}{\kappa_{0}+j+\frac{1}{2},2j+1}%
;1\right)  .
\]
Then%
\begin{gather*}
g_{1}\left(  z\right)     =\sum_{j=-n}^{n}\left(  n-m+1\right)  _{\left\vert
j\right\vert -1}\frac{n!}{\left(  n-\left\vert j\right\vert \right)
!}b_{\left\vert j\right\vert }\left(  n-m+j\right)  z^{n+j}\overline{z}%
^{n-j},\\
g_{2}\left(  z\right)     =\sum_{j=n-m}^{m-n}\left(  -n\right)  _{\left\vert
j\right\vert }\frac{\left(  m-n-1\right)  !}{\left(  m-n-\left\vert
j\right\vert \right)  !}b_{\left\vert j\right\vert }\left(  m-n+j\right)
\overline{z}^{n+j}z^{n-j}.
\end{gather*}
Thus $c\left(  g_{1};n,n\right)  =b_{0}=c\left(  g_{2};n,n\right)  $. Suppose
$\left\vert j\right\vert \geq1$, then
\begin{gather*}
c\left(  g_{1};n+j,n-j\right)  =
\left(
n-m+1\right)  _{\left\vert j\right\vert -1}\left(  -n\right)  _{\left\vert
j\right\vert }\left(  -1\right)  ^{j}\left(  n-m+j\right)  b_{\left\vert
j\right\vert }
\end{gather*}
for $\left\vert j\right\vert \leq n$, and the equation remains
valid if $n<\left\vert j\right\vert \leq m-n$ because $\left(  -n\right)
_{\left\vert j\right\vert }=0$ for $\left\vert j\right\vert >n$. Also
$c\left(  g_{2};n+j,n-j\right)  =\left(  -n\right)  _{\left\vert j\right\vert
}\left(  n+1-m\right)  _{\left\vert j\right\vert -1}\left(  -1\right)
^{j-1}\left(  m-n-j\right)  b_{\left\vert j\right\vert }$, and the equation
remains valid if $m-n<\left\vert j\right\vert \leq n$ (that is $n-m+\left\vert
j\right\vert -1\geq0$). Thus $c\left(  g_{1};n+j,n-j\right)  =c\left(
g_{2};n+j,n-j\right)  $ for $\left\vert j\right\vert \leq\max\left(
n,m-n\right)  $ and $g_{1}=g_{2}$.
\end{proof}

Note that if $k=0,1,2,\ldots$ then $\left(  -k\right)  _{j}=0$ for $j>k$.
Recall the structural constants for the Poisson kernel $P_{n}$ from equations
$($\ref{lb02n})--(\ref{lb2n1}). These are rational functions of $\kappa
_{0}$, $\kappa_{1}$ def\/ined for all $\left(  \kappa_{0},\kappa_{1}\right)  \notin
K_{0}$.

\begin{proposition}
\label{lbm2n}Suppose $-\left(  \kappa_{0}+\kappa_{1}\right)  =m\in%
\mathbb{N}
$, $1\leq n\leq m-1,$ and $n\neq\frac{m}{2}$ then%
\begin{gather*}
\frac{\lambda_{2m-2n}^{0}}{\lambda_{2n}^{0}}    =-\left(  \frac{\left(
\kappa_{0}+\frac{1}{2}\right)  _{n}\left(  m-n\right)  !}{\left(  \kappa
_{0}+\frac{1}{2}\right)  _{m-n}n!}\right)  ^{2},\qquad
\frac{\lambda_{2m-2n}^{1}}{\lambda_{2n}^{1}}    =-\left(  \frac{\left(
\kappa_{0}+\frac{1}{2}\right)  _{n}\left(  m-n-1\right)  !}{\left(  \kappa
_{0}+\frac{1}{2}\right)  _{m-n}\left(  n-1\right)  !}\right)  ^{2},\\
\lambda_{2m}^{0}    =-\left(  \frac{m!}{\left(  \kappa_{0}+\frac{1}%
{2}\right)  _{m}}\right)  ^{2},\qquad \lambda_{2m}^{1}=0.
\end{gather*}
\end{proposition}

\begin{proof}
Recall $\lambda_{2n}^{0}=\frac{n!\left(  \kappa_{0}+\kappa_{1}+1\right)
_{n-1}\left(  \kappa_{0}+\kappa_{1}+2n\right)  }{\left(  \kappa_{0}+\frac
{1}{2}\right)  _{n}\left(  \kappa_{1}+\frac{1}{2}\right)  _{n}}$ (for $n\in%
\mathbb{N}
$). Also $\left(  \kappa_{1}+\frac{1}{2}\right)  _{n}=\allowbreak\left(
-m-\kappa_{0}+\frac{1}{2}\right)  _{n}=\allowbreak\left(  -1\right)
^{n}\left(  \kappa_{0}+\frac{1}{2}+m-n\right)  _{n}=\allowbreak\left(
-1\right)  ^{n}\frac{\left(  \kappa_{0}+\frac{1}{2}\right)  _{m}}{\left(
\kappa_{0}+\frac{1}{2}\right)  _{m-n}}$, and similarly $\left(  \kappa
_{1}+\frac{1}{2}\right)  _{m-n}=\left(  -1\right)  ^{m-n}\frac{\left(
\kappa_{0}+\frac{1}{2}\right)  _{m}}{\left(  \kappa_{0}+\frac{1}{2}\right)
_{n}}$. Thus%
\[
\frac{\lambda_{2m-2n}^{0}}{\lambda_{2n}^{0}}=\left(  -1\right)  ^{m}\left(
\frac{\left(  \kappa_{0}+\frac{1}{2}\right)  _{n}}{\left(  \kappa_{0}+\frac
{1}{2}\right)  _{m-n}}\right)  ^{2}\frac{\left(  m-n\right)  !\left(
1-m\right)  _{m-n-1}\left(  m-2n\right)  }{n!\left(  1-m\right)  _{n-1}\left(
-m+2n\right)  }.
\]
But $\frac{\left(  1-m\right)  _{m-n-1}}{\left(  1-m\right)  _{n-1}}=\left(
-1\right)  ^{m}\frac{\left(  m-1\right)  !\left(  m-n\right)  !}{\left(
m-1\right)  !n!}$ (note $\left(  -k\right)  _{j}=\left(  -1\right)  ^{j}%
\frac{k!}{\left(  k-j\right)  !}$ for $k\in%
\mathbb{N}
_{0}$). Next $\lambda_{2n}^{1}=\frac{\left(  n-1\right)  !\left(  \kappa
_{0}+\kappa_{1}+1\right)  _{n}\left(  \kappa_{0}+\kappa_{1}+2n\right)
}{\left(  \kappa_{0}+\frac{1}{2}\right)  _{n}\left(  \kappa_{1}+\frac{1}%
{2}\right)  _{n}}$. Similarly we f\/ind%
\[
\frac{\lambda_{2m-2n}^{1}}{\lambda_{2n}^{1}}=\left(  -1\right)  ^{m}\left(
\frac{\left(  \kappa_{0}+\frac{1}{2}\right)  _{n}}{\left(  \kappa_{0}+\frac
{1}{2}\right)  _{m-n}}\right)  ^{2}\frac{\left(  m-n-1\right)  !\left(
1-m\right)  _{m-n}\left(  m-2n\right)  }{\left(  n-1\right)  !\left(
1-m\right)  _{n}\left(  -m+2n\right)  },
\]
and $\frac{\left(  1-m\right)  _{m-n}}{\left(  1-m\right)  _{n}}=\left(
-1\right)  ^{m}\frac{\left(  m-1\right)  !\left(  m-n-1\right)  !}{\left(
m-1\right)  !\left(  n-1\right)  !}$. The special case $\lambda_{2m}^{0}$
follows from setting $n=0$ in the f\/irst formula. The term $\left(  \kappa
_{0}+\kappa_{1}+1\right)  _{m}$ shows $\lambda_{2m}^{1}=0$.
\end{proof}

The following two propositions are proven by similar calculations.

\begin{proposition}
\label{lbm2n1}Suppose $-\left(  \kappa_{0}+\kappa_{1}\right)  =m\in%
\mathbb{N}
$, $0\leq n\leq m-1,$ and $n\neq\frac{m-1}{2}$ then%
\begin{gather*}
\frac{\lambda_{2m-2n-1}^{0}}{\lambda_{2n+1}^{0}}    =-\left(  \frac{\left(
\kappa_{0}+\frac{1}{2}\right)  _{n+1}\left(  m-n-1\right)  !}{\left(
\kappa_{0}+\frac{1}{2}\right)  _{m-n}n!}\right)  ^{2},
\\
\frac{\lambda_{2m-2n-1}^{1}}{\lambda_{2n+1}^{1}}    =-\left(  \frac{\left(
\kappa_{0}+\frac{1}{2}\right)  _{n}\left(  m-n-1\right)  !}{\left(  \kappa
_{0}+\frac{1}{2}\right)  _{m-n-1}n!}\right)  ^{2}.
\end{gather*}
\end{proposition}

\begin{proposition}
\label{lbm00}Suppose $-\left(  \kappa_{0}+\kappa_{1}\right)  =m\in%
\mathbb{N}
$, $0\leq n\leq m-1$ then%
\[
\frac{\lambda_{2m-2n-1}}{\lambda_{2n}}=-\left(  \frac{\left(  \kappa_{0}%
+\frac{1}{2}\right)  _{n}\left(  m-n-1\right)  !}{\left(  \kappa_{0}+\frac
{1}{2}\right)  _{m-n}n!}\right)  ^{2}.
\]

\end{proposition}

\begin{proposition}
\label{lbzero}Suppose $-\left(  \kappa_{0}+\kappa_{1}\right)  =m\in%
\mathbb{N}
$ then $\lambda_{m}^{0}=0=\lambda_{m}^{1}$. If $n>2m$ then $\lambda_{n}%
^{0}=0=\lambda_{n}^{1}$. If $n\geq2m$ then $\lambda_{n}=0$.
\end{proposition}

\begin{proof}
Since both $\lambda_{n}^{0}$ and $\lambda_{n}^{1}$ contain the factor $\left(
\kappa_{0}+\kappa_{1}+n\right)  $ for $n$ even or odd, it follows that
$\lambda_{m}^{0}=0=\lambda_{m}^{1}$. The term $\left(  \kappa_{0}+\kappa
_{1}+1\right)  _{j}$ vanishes for $j>m$, and $j=k-1$ for $\lambda_{2k}^{0}$,
$j=k$ for each of $\lambda_{2k}$, $\lambda_{2k+1}$, $\lambda_{2k}^{1}$,
$\lambda_{2k+1}^{0}$, $\lambda_{2k+1}^{1}$.
\end{proof}

\begin{theorem}
\label{Pzero}Suppose $-\left(  \kappa_{0}+\kappa_{1}\right)  =m\in%
\mathbb{N}
$ then $P_{N}\left(  z,w\right)  =0$ for $N>2sm$ and $P_{N}\left(  z,w\right)
+\left(  z\overline{z}w\overline{w}\right)  ^{N-sm}P_{2sm-N}\left(
z,w\right)  =0$ for $0\leq N\leq2sm$.
\end{theorem}

\begin{proof}
If $N=sk>2sm$ then $P_{sk}\left(  z,w\right)  =\lambda_{k}^{0}f_{k}^{0}\left(
z^{s}\right)  f_{k}^{0}\left(  w^{s}\right)  +\lambda_{k}^{1}f_{k}^{1}\left(
z^{s}\right)  f_{k}^{1}\left(  \overline{w}^{s}\right)  $ and $\lambda_{k}%
^{0}=\lambda_{k}^{1}=0$ by Proposition \ref{lbzero}. If $N=sk+t$ with $1\leq
t<s$ and $N>2sm$ then $P_{sk+t}\left(  z,w\right)  =\lambda_{k}\left(
z^{t}\overline{w}^{t}f_{k}\left(  z^{s}\right)  f_{k}\left(  \overline{w}%
^{s}\right)  +\overline{z}^{t}w^{t}f_{k}\left(  \overline{z}^{s}\right)
f_{k}\left(  w^{s}\right)  \right)  $, $k\geq2sm$ and $\lambda_{k}=0$. If
$N=sm$ then $\lambda_{m}^{0}=\lambda_{m}^{1}=0$. Suppose $N=2ns$ and $0<N<2sm$
(so $0<n<m$), then%
\begin{gather*}
P_{2ns}\left(  z,w\right)     =\lambda_{2n}^{0}\left(  \frac{\left(
\kappa_{0}+\frac{1}{2}\right)  _{n}\left(  m-n\right)  !}{\left(  \kappa
_{0}+\frac{1}{2}\right)  _{m-n}n!}\right)  ^{2}\left(  z\overline{z}%
w\overline{w}\right)  ^{\left(  2n-m\right)  s}f_{2m-2n}^{0}\left(
z^{s}\right)  f_{2m-2n}^{0}\left(  w^{s}\right) \\
\phantom{P_{2ns}\left(  z,w\right)     =}{}  +\lambda_{2n}^{1}\left(  \frac{\left(  \kappa_{0}+\frac{1}{2}\right)
_{n}\left(  m-n-1\right)  !}{\left(  \kappa_{0}+\frac{1}{2}\right)
_{m-n}\left(  n-1\right)  !}\right)  ^{2}\left(  z\overline{z}w\overline
{w}\right)  ^{\left(  2n-m\right)  s}f_{2m-2n}^{1}\left(  z^{s}\right)
f_{2m-2n}^{1}\left(  \overline{w}^{s}\right)  ,
\end{gather*}
thus by\ Propositions \ref{f2n0m} and \ref{lbm2n}
\begin{gather*}
\left(  z\overline{z}w\overline{w}\right)  ^{\left(  m-2n\right)  s}%
P_{2ns}\left(  z,w\right)  +P_{2ms-2ns}\left(  z,w\right)  \\
\qquad{}=\lambda_{2n}^{0}\left\{  \left(  \frac{\left(  \kappa_{0}+\frac{1}{2}\right)
_{n}\left(  m-n\right)  !}{\left(  \kappa_{0}+\frac{1}{2}\right)  _{m-n}%
n!}\right)  ^{2}+\frac{\lambda_{2m-2n}^{0}}{\lambda_{2n}^{0}}\right\}
f_{2m-2n}^{0}\left(  z^{s}\right)  f_{2m-2n}^{0}\left(  w^{s}\right) \\
\qquad{}+\lambda_{2n}^{1}\left\{  \left(  \frac{\left(  \kappa_{0}+\frac{1}{2}\right)
_{n}\left(  m-n-1\right)  !}{\left(  \kappa_{0}+\frac{1}{2}\right)
_{m-n}\left(  n-1\right)  !}\right)  ^{2}+\frac{\lambda_{2m-2n}^{1}}%
{\lambda_{2n}^{1}}\right\}  f_{2m-2n}^{1}\left(  z^{s}\right)  f_{2m-2n}%
^{1}\left(  \overline{w}^{s}\right) =0.
\end{gather*}
For the special case $N=2sm$ we have $P_{2sm}\left(  z,w\right)  =\lambda
_{2m}^{0}\Big(  \frac{\left(  \kappa_{0}+\frac{1}{2}\right)  _{m}}%
{m!}\Big)  ^{2} \!\! \left(  z\overline{z}w\overline{w}\right)  ^{ms}%
\allowbreak=-\left(  z\overline{z}w\overline{w}\right)  ^{ms}P_{0}$ because
$\lambda_{2m}^{1}=0,P_{0}=1$ and $\lambda_{2m}^{0}=-\left(  \frac{m!}{\left(
\kappa_{0}+\frac{1}{2}\right)  _{m}}\right)  ^{2}$. Similarly\ by use of
Propositions \ref{f2n1m} and \ref{lbm2n1} we show the result holds for
$N=\left(  2n+1\right)  s$ for $0\leq n<m$ and $2n+1\neq m$. Suppose $N=sk+t$
with $1\leq t<s$ and $0<N<2sm$, then $2sm-N=s\left(  2m-k-1\right)  +\left(
s-t\right)  $. One of $k,2m-k-1$ is even so assume $k=2n$ with $0\leq n<m$
(otherwise replace $N$ by $2sm-N$ and $t$ by $s-t$). By Propositions
\ref{f2nm1} and \ref{lbm00}%
\begin{gather*}
  \lambda_{2n}z^{t}\overline{w}^{t}f_{2n}\left(  z^{s}\right)  f_{2n}\left(
\overline{w}^{s}\right)
  {+}\lambda_{2m-2n-1}\left(  z\overline{z}w\overline{w}\right)  ^{\left(
2n-m\right)  s+t}\overline{z}^{s-t}w^{s-t}f_{2m-2n-1}\left(  \overline{z}%
^{s}\right)  f_{2m-2n-1}\left(  w^{s}\right)\! \\
\qquad{}  =\lambda_{2n}\left(  z\overline{w}\right)  ^{\left(  2n-m\right)
s+t}\left(  \overline{z}w\right)  ^{\left(  2n-m+1\right)  s}f_{2m-2n-1}%
\left(  \overline{z}^{s}\right)  f_{2m-2n-1}\left(  w^{s}\right) \\
\qquad{} \times\left\{  \left(  \frac{\left(  \kappa_{0}+\frac{1}{2}\right)
_{n}\left(  m-n-1\right)  !}{\left(  \kappa_{0}+\frac{1}{2}\right)  _{m-n}%
n!}\right)  ^{2}+\frac{\lambda_{2m-2n-1}}{\lambda_{2n}}\right\}   =0.
\end{gather*}
Add this equation to its complex conjugate to show
\begin{gather*}
P_{2ns+t}\left(
z,w\right)  +\left(  z\overline{z}w\overline{w}\right)  ^{\left(  2n-m\right)
s+t}P_{\left(  2m-2n\right)  s-t}\left(  z,w\right)  =0.\tag*{\qed}
\end{gather*}\renewcommand{\qed}{}
\end{proof}

\begin{theorem}
For $n,m{\in}
\mathbb{N}
$ equation \eqref{Kser} for $K_{n}\!\left(  z,w\right)  $ has a removable
singularity at $\kappa_{0}{+}\kappa_{1}{=}{-}m.\!$
\end{theorem}

\begin{proof}
Consider the series
\[
K_{n}\left(  z,w\right)  =2^{-n}\sum_{j=0}^{\left\lfloor n/2\right\rfloor
}\frac{1}{j!\left(  s\kappa_{0}+s\kappa_{1}+1\right)  _{n-j}}\left(
z\overline{z}w\overline{w}\right)  ^{j}P_{n-2j}\left(  z,w\right)  .
\]
The possible poles occur at $n-j\geq sm$ (that is, $\left(  1-sm\right)
_{n-j}=0$) and the multiplicities do not exceed $1$. Thus there are no poles
if $n<sm$. If $n-2j>2sm$ then $P_{n-2j}\left(  z,w\right)  $ is divisible by
$\left(  \kappa_{0}+\kappa_{1}+m\right)  $, by Theorem \ref{Pzero}, and the
singularity is removable. It remains to consider the case $n-2j\leq2sm$ and
$n-j\geq sm$. Suppose $j=j_{0}$ satisf\/ies these inequalities and let
$j_{1}=n-j_{0}-sm$. Then $j_{1}\geq0$ and $n-2j_{1}=2sm-n+2j_{0}\geq0$, hence
$j=j_{1}$ appears in the sum. But $2sm-\left(  n-2j_{0}\right)  =n-2j_{1}$ so
Theorem \ref{Pzero} applies. We can assume $j_{1}\leq j_{0}$. Consider the
following subset of the sum for $K_{n}\left(  z,w\right)  $:%
\begin{gather*}
  \frac{\left(  z\overline{z}w\overline{w}\right)  ^{j_{0}}P_{n-2j_{0}%
}\left(  z,w\right)  }{j_{0}!\left(  s\kappa_{0}+s\kappa_{1}+1\right)
_{n-j_{0}}}+\frac{\left(  z\overline{z}w\overline{w}\right)  ^{j_{1}%
}P_{n-2j_{1}}\left(  z,w\right)  }{j_{1}!\left(  s\kappa_{0}+s\kappa
_{1}+1\right)  _{n-j_{1}}}
  =\frac{\left(  z\overline{z}w\overline{w}\right)  ^{j_{0}}}{j_{0}!\left(
s\kappa_{0}+s\kappa_{1}+1\right)  _{n-j_{0}}}C_{n,j_{0}}.
\end{gather*}
with
\[
C_{n,j_{0}}=P_{n-2j_{0}}\left(  z,w\right)  +\frac{j_{0}!\left(  z\overline
{z}w\overline{w}\right)  ^{j_{1}-j_{0}}P_{n-2j_{1}}\left(  z,w\right)  }%
{j_{1}!\left(  s\kappa_{0}+s\kappa_{1}+1+n-j_{0}\right)  _{j_{0}-j_{1}}}.
\]
The expression $C_{n,j_{0}}$ has no pole at $\kappa_{0}+\kappa_{1}=-m$ since
$1-sm+n-j_{0}\geq1$. Indeed $\left(  1-sm+n-j_{0}\right)  _{j_{0}-j_{1}%
}=\left(  j_{1}+1\right)  _{j_{0}-j_{1}}=j_{0}!/j_{1}!$. In the special case
$n-2j_{0}=ms$, and $j_{0}=j_{1}=\left(  n-sm\right)  /2$ we replace
$C_{n,j_{0}}$ by $P_{n-2j_{0}}\left(  z,w\right)  $. By Theorem \ref{Pzero}
$C_{n,j_{0}}=0$ when $\kappa_{0}+\kappa_{1}=-m$, thus $C_{n,j_{0}}$ is
divisible by $\left(  \kappa_{0}+\kappa_{1}+m\right)  s$. The sum of the two
terms ($j=j_{0}$ and $j=j_{1}$) has a removable singularity there.
\end{proof}

The expressions for $V\left(  z^{a}\overline{z}^{b}\right)  $ are derived from
the series \eqref{Kser} for $K_{n}\left(  z,w\right)  $ thus the result about
singularities at $\kappa_{0}+\kappa_{1}=-m$ being removable by grouping the
expansion into certain pairs applies. Note that in the above proof the paired
terms are $P_{n-2j_{0}}$ and $P_{n-2j_{1}}$ with $j_{0}+j_{1}=n-sm$. To
analyze $V\left(  z^{a}\overline{z}^{b}\right)  $ it suf\/f\/ices to identify the
pairs. For the case $a\equiv b\operatorname{mod}s$ and $a\geq b$ let
$a=us+r$, $b=vs+r$ and $0\leq r<s$. The paired indices in the sum from Theorem
\ref{Vz2} consist of $\left\{  \left(  k,k^{\prime}\right)  :0\leq
k<k^{\prime}\leq v,k+k^{\prime}=u+v-m\right\}  $ . Indeed for $k$, $k^{\prime}$
with $0\leq k,k^{\prime}\leq v$ def\/ine $j$ by $a+b-2j=\left(  u+v-2k\right)
s$ so that $j=ks+r$ and similarly set $j^{\prime}:=k^{\prime}s+r$. The pairing
condition $j+j^{\prime}=a+b-sm$ is equivalent to $k+k^{\prime}=u+v-m$. Thus
$k$, $k^{\prime}$ are paired exactly when $k+k^{\prime}=u+v-m$ and $0\leq
k,k^{\prime}\leq v$.

For the case $a-b\equiv t\operatorname{mod}s$, and with $a=us+r+t>b=vs+r$,
$0\leq r<s$, $1\leq t<s$ the pairing in the formula from Theorem \ref{Vz1}
combines terms from the f\/irst sum with corresponding terms in the second. For
the f\/irst sum suppose $0\leq k\leq v$ and $j:=ks+r$ so that $a+b-2j=t+\left(
u+v-2k\right)  s$. For the second sum let $1-\left\lfloor \frac{r+t}%
{s}\right\rfloor \leq k^{\prime}\leq v$ and let $j^{\prime}:=\left(
k^{\prime}-1\right)  s+r+t$ so that $a+b-2j^{\prime}=\left(  s-t\right)
+\left(  u+v-2k^{\prime}+1\right)  s$. The pairing condition $j+j^{\prime
}=a+b-sm$ is equivalent to $k+k^{\prime}=u+v+1-m$. To remove the singularities
at $\kappa_{0}+\kappa_{1}=-m$ combine the term in the f\/irst sum of index $k$
with the term in the second of index $k^{\prime}$ for all pairs $\left(
k,k^{\prime}\right)  $ satisfying $k+k^{\prime}=u+v+1-m$, $0\leq k\leq v$,
$1-\left\lfloor \frac{r+t}{s}\right\rfloor \leq k^{\prime}\leq v$.

\pdfbookmark[1]{References}{ref}
\LastPageEnding

\end{document}